\documentclass[12pt]{elsarticle}
\usepackage{amsmath,amssymb,amsthm,amsfonts}
\usepackage[utf8]{inputenc}
\usepackage[english]{babel}
\usepackage{geometry}
\usepackage{graphicx}
\usepackage{multirow}
\usepackage{tikz}
\usetikzlibrary{shapes,arrows,positioning}
\usepackage{appendix}
\usepackage{lscape}
\usepackage{hyperref}
\usepackage{bbm}

\theoremstyle{plain}

\theoremstyle{definition}

\theoremstyle{remark}

\begin{document}

\begin{frontmatter}
\title{Comparative study on higher order compact RBF-FD formulas with Gaussian and Multiquadric radial functions}

\author[inst1]{Manoj Kumar Yadav\corref{cor1}}
\cortext[cor1]{Corresponding author}
\ead{manojkumar.yadav@mahindrauniversity.edu.in}

\author[inst1]{Chirala Satyanarayana}
\ead{satyanarayana.chirala@mahindrauniversity.edu.in}

\author[inst1]{A. Sreedhar}
\ead{sreedhar22pmat001@mahindrauniversity.edu.in}

\affiliation[inst1]{organization={Department of Mathematics}, {Ecole Centrale School of Engineering, Mahindra University}, 
           city={Hyderabad},
           postcode={500043},
           state={Telangana},
           country={India.}}
\begin{abstract}
\noindent 
We generate Gaussian radial function based higher order compact RBF-FD formulas for some differential operators. Analytical expressions for weights associated to first and second derivative formulas (up to order 10) and 2D-Laplacian formulas (up to order 6) are derived. Then these weights are used to obtain analytical expression for local truncation errors. The weights are obtained by symbolic computation of a linear system in {\it Mathematica}. Often such linear systems are not directly amenable to symbolic computation. We make use of symmetry of formula stencil along with Taylor series expansions for performing the computation. In the flat limit, the formulas converge to their respective order polynomial based compact FD formulas. We validate the formulas with standard test functions and demonstrate improvement in approximation accuracy with respect to corresponding order multiquadric based compact RBF-FD formulas and compact FD schemes. We also compute optimal value of shape parameter for each formula.
\end{abstract}
\end{frontmatter}
\section{Introduction}\label{sec1}
To approximate differential operators, Wright and Fornberg \cite{wright2006scattered} proposed Radial Basis Function based Hermite Finite Difference (RBF-HFD) formulas -- also referred as compact RBF-FD formulas -- comprising linear combination of the unknown function and its operator values evaluated on unstructured nodes. They analyzed the behaviour of Multiquadric (MQ) based RBF-FD and RBF-HFD formulas in the flat limit and observed convergence to FD and compact FD formulas reported in Collatz \cite{collatznumerical} and Lele \cite{lele1992compact}. They also observed ill-conditioning in the numerical computations due to RBF-Direct approach and used Contour-Padé algorithm \cite{fornberg2004stable} for stable computations in the flat limit. Later, several other stable algorithms \cite{fornberg2011stable, fornberg2013stable, wright2017stable} were developed for flat limit computations. 
To address the problem of ill-conditioning, Bayona {\it et al.} \cite{bayona2010rbf} obtained analytical expressions for MQ based RBF-FD formulas and their convergence to classical FD formulas. Bayona and coworkers \cite{bayona2011optimal, bayona2012optimal} applied the formulas derived in \cite{bayona2010rbf} to solve PDE problems with constant and variable shape parameters. Bayona, Moscoso, and Kindelan \cite{bayona2012gaussian} derived analytical expressions for weights of Gaussian (GA) RBF-FD formulas for first derivative, second derivative and 2D-Laplacian. For first and second derivative they obtained fourth order GA based RBF-HFD formulas. \\ ~ \\
\noindent Satyanarayana {\it et al.} \cite{Satya_RBF_HFD} derived leading/higher order analytical expressions for weights and local truncation errors (LTEs) of MQ based RBF-HFD formulas for first and second derivatives (up to tenth order) and for 2D-Laplacian up to sixth order. The formulas were validated for some test functions and convergence to respective order FD schemes were established. The analytical expressions were helpful in removing ill-conditioning from the approximations and locating optimal value of shape parameter for a chosen test function. More recently, Song {\it et al.} \cite{song2024computing}, made use of local collocation on non-uniform nodes set with RBFs and integrated RBFs augmented with polynomials to derive analytical expressions for compact FD type formulas for first derivative, second derivative and 2D-Laplacian.\\~ \\
\noindent  In this paper, we extend the work of Bayona {\it et al.} \cite{bayona2012gaussian} and Satyanarayana {\it et al.} \cite{Satya_RBF_HFD} to derive analytical expressions for weights of Gaussian based compact FD (RBF-HFD) formulas for first and second derivative up to tenth order and up to sixth order approximation for 2D-Laplacian. We establish convergence of the approximation formulas to the respective polynomial based compact FD formulas in the limit of $\epsilon \rightarrow 0 $ by comparing the limiting weights with those reported in Collatz \cite{collatznumerical} and Lele \cite{lele1992compact}. Further, we compare the approximation accuracy of GA based formulas with Multiquadric (MQ) based formulas and compact FD schemes \cite{collatznumerical, lele1992compact} for various test functions \cite{wright2006scattered, bayona2010rbf, ding2005error, chandhini2007local}. Finally, we compute optimal values of shape parameter for each formula following an optimization technique originally proposed by Bayona {\it et al.} \cite{bayona2010rbf} for RBF-FD approximations. \\~ \\
\noindent Organization of the article is as follows. In Sections 2 and 3, we derive Gaussian based RBF-HFD formulas (order 6, 8, 10) for first and second derivative. Analytical expressions for weights and local truncation errors associated to the formulas are obtained. We also validate The formulas for standard test functions and establish their convergence properties in the limit of $\epsilon$ tending to zero. Section 4, presents analytical expressions of weights and local truncation errors associated to fourth and sixth order GA based RBF-HFD formulas to approximate 2D Laplacian. In Section 5, we compute optimal values of shape parameter for the formulas obtained in the previous sections. Section 6 presents conclusions of the present work.
\section{RBF-HFD method} \label{sec2_RBF_HFD_method}
                
        \noindent For the sake of completeness we briefly describe RBF-HFD method. One may refer Satyanarayana {\it et al.} \cite{Satya_RBF_HFD} for more detailed description. Let $\phi(\mathbf{x},\mathbf{y}):=\phi(\|\mathbf{x}-\mathbf{y}\|_2)$ denote a radial function centered at node $\mathbf{y}$. In this paper, we perform all the derivations {\it w.r.t.} Gaussian (GA) radial function. Therefore, going forward $\phi(\|\mathbf{x}-\mathbf{y}\|_2)=e^{-\epsilon^2 \|\mathbf{x}-\mathbf{y}\|_2^2}$, where $\epsilon>0$ is a shape parameter which determines the curvature of the radial function. For a fixed reference node $\mathbf{x}_0 \in \mathbb{R}^{d}$, let $S=\{\mathbf{x}_{1},\mathbf{x}_{2},\ldots,\mathbf{x}_{n} \}$ be a local support and $\mu = \{ \mathbf{x}_{1},\ldots,\mathbf{x}_{-1},\mathbf{x}_{1},\ldots,\mathbf{x}_{m} \} $ with $m \leq n$ denote deleted local support of $\mathbf{x}_0$. Then in view of RBF-HFD method \cite{wright2006scattered}, a linear differential operator ($\mathcal{L}u$) is approximated at $\mathbf{x}_0$ by linear combination of function values $u$ over set $S$ and its operator values $\mathcal{L}u$ over set $\mu$. Thus, the RBF- HFD approximation formula may be written as 
 \begin{equation}
            S(x_{i})=u(x_{i})~ ~\text{for} ~ ~i= 1\cdots N, ~ ~ 
            \mathcal{L}S(x_{k})=\mathcal{L}u(x_{k}) ~ ~\text{for} ~ ~k= 1\cdots m.
        \end{equation}  
        where interpolation takes the form as
        \begin{equation}\label{RBF_HFD_Method}
           S(x_{i}) \approx \sum_{i=1}^{n} \alpha_i u_i + \sum_{k=1\;(k\;\neq \; r)}^{m} \beta_{k}\mathcal{L}u_k+ \sum_{l=1}^{M} d_{l}P_{l}(x)~ ~ 
           \text{where} ~ ~ u_i=u(\mathbf{x}_i),~ ~ \mathcal{L}u_k=\mathcal{L}u(\mathbf{x}_k).
        \end{equation}
         \noindent The unknown weights $\alpha:=\{\alpha_i:1\leq i\leq n\}$ and $\beta:=\{\beta_k:1\leq k \leq m,k\neq 0\}$ in formula (\ref{RBF_HFD_Method}) are computed by collocating the formula on sets of radial basis functions
        $$B_S:=\{\phi_j(\mathbf{x})=\phi(\|\mathbf{x}-\mathbf{x}_j\|_2): 1\leq j \leq n\}~~\mbox{and}~~B_{\mu}:=\{\mathcal{L}^2\phi_k(\mathbf{x}): 1\leq k\leq m, k\neq 0\}.$$ 
          The collocation procedure result in the following linear system
\begin{equation} \label{RBF_HFD_Sys}
         \begin{bmatrix}
		  \Phi & {\mathcal{L}^2}\Phi & \textbf{P} \\
            {\mathcal{L}^1}\Phi  & {{\mathcal{L}^1}{\mathcal{L}^2}}\Phi & {\mathcal{L}^1\textbf{P}} \\
            \textbf{P}^T & {\mathcal{L}^1}\textbf{P}^T & 0
\end{bmatrix} 
\begin{bmatrix}
    \alpha \\
    \beta \\
    \gamma
\end{bmatrix} = \begin{bmatrix}
    u(x_{i}) \\
    \mathcal{L}^1 u(x_{i}) \\
    0
\end{bmatrix}.
\end{equation}
The linear system of equations (\ref{RBF_HFD_Sys}) is solved using symbolic computation in {\it Mathematica} to obtain analytical expressions for the weights $\alpha$, $\beta$. These weights are then substituted in equation (\ref{RBF_HFD_LTE}) and expanded using Taylor series expansion in powers of $\epsilon h$ for $\epsilon h \ll 1$ to obtain leading order expression for the local truncation error ($\tau_0$) 
\begin{equation}\label{RBF_HFD_LTE}
    \tau_{0} =\left( \sum_{i=1}^{n} \alpha_i u_i + \sum_{k=1 \;(k \; \neq \; 0) }^{m} \beta_{k} \mathcal{L}u_k +\sum_{l=1}^{M} d_{l}P_{l}(x)\right) -\mathcal{L}u(\mathbf{x}_r). 
\end{equation} 
The block matrix entries in the linear system (\ref{RBF_HFD_Sys}) are defined as \\
\begin{equation*}   
\phi_{N \times N}=\begin{pmatrix}
     \phi(\|x_{1}-x_{1}\|) & \phi(\|x_{1}-x_{2}\|) & \phi(\|x_{1}-x_{3}\|) & \cdots & \phi(\|x_{1}-x_{n}\|) \\
    \phi(\|x_{2}-x_{1}\|) & \phi(\|x_{2}-x_{2}\|) & \phi(\|x_{2}-x_{3}\|) & \cdots & \phi(\|x_{2}-x_{n}\|) \\
    \vdots & \vdots & \vdots &\cdots & \vdots \\
    \phi(\|x_{n}-x_{1}\|) & \phi(\|x_{n}-x_{2}\|) & \phi(\|x_{n}-x_{3}\|) & \cdots & \phi(\|x_{n}-x_{n}\|)
\end{pmatrix}
\end{equation*}
\vspace{0.2cm}
\begin{equation*}   
{\mathcal{L}^2}\Phi_{N \times M}=\begin{pmatrix}
    {\mathcal{L}^2}\phi(\|x_{1}-x_{1}\|) & {\mathcal{L}^2}\phi(\|x_{1}-x_{2}\|) & {\mathcal{L}^2}\phi(\|x_{1}-x_{3}\|) & \cdots & {\mathcal{L}^2}\phi(\|x_{1}-x_{n}\|) \\
    {\mathcal{L}^2}\phi(\|x_{2}-x_{1}\|) & {\mathcal{L}^2}\phi(\|x_{2}-x_{2}\|) & {\mathcal{L}^2}\phi(\|x_{2}-x_{3}\|) & \cdots & {\mathcal{L}^2}\phi(\|x_{2}-x_{n}\|)\\
    \vdots & \vdots & \vdots & \cdots & \vdots \\
    {\mathcal{L}^2}\phi(\|x_{n}-x_{1}\|) & {\mathcal{L}^2}\phi(\|x_{n}-x_{2}\|) & {\mathcal{L}^2}\phi(\|x_{n}-x_{3}\|) & \cdots & {\mathcal{L}^2}\phi(\|x_{n}-x_{n}\|)
\end{pmatrix}
\end{equation*}
\vspace{0.2cm}
\begin{equation*}   
P_{N \times (N-2)}=\begin{pmatrix}
    1 & x_{1} & y_{1} \\
    1 & x_{2} & y_{2} \\
    \vdots & \vdots & \vdots \\
    1 & x_{n} & y_{n} \\
\end{pmatrix}, ~ ~
{\mathcal{L}^1}P_{M \times (N-2)}=\begin{pmatrix}
    0 & {\mathcal{L}^1}(x_{1}) &{\mathcal{L}^1}(y_{1}) \\
    0 & {\mathcal{L}^1}(x_{2}) &{\mathcal{L}^1}(y_{2}) \\
    \vdots & \vdots & \vdots \\
    0 & {\mathcal{L}^1}(x_{n}) &{\mathcal{L}^1}(y_{n}) \\
\end{pmatrix}
\end{equation*}
\section{First Derivative Approximations} \label{First_derivative_approximation}
\noindent In this section, we obtain analytical expressions of weights and local truncation errors for sixth, eighth and tenth order GA based RBF-HFD formulas for first derivative. We consider uniform distribution of nodes in the formula stencils. Therefore, the formulas depend on stencil parameter $h$ and shape parameter $\epsilon$ associated with the Gaussian RBF. Symmetry of uniform stencil and central difference approximation of first derivative result in fewer unknown weights in the formulas.
The obtained formulas converge to respective order compact FD formulas \cite{collatznumerical,lele1992compact} in the limit of $\epsilon$ tending to zero. We validate the formulas' approximation efficiency {\it w.r.t.} the following test functions considered in Bayona {\it et al.} \cite{bayona2012gaussian} and Chandhini {\it et al.} \cite{sanyasiraju2008local}  
\begin{eqnarray}
&&u_{1}(x)=\sin({x}^2),  \;\;\; x_0=0.4, \label{u1} \\
&&u_{2}(x)=\sin(\pi x ) + \dfrac{e^{x}-1}{e-1},  \;\;\; x_0=0.25. \label{u2}
\end{eqnarray}
Further, we compare the approximation accuracy of GA based formulas with MQ based formulas \cite{Satya_RBF_HFD} and compact FD schemes \cite{collatznumerical, lele1992compact}.
\subsection{Sixth order formula}
\noindent The sixth order formula stencil has local support sets  $S=\{  x_0-2h, x_0-h, x_0, x_0+h, x_0+2h  \} $ and $\mu=\{x_0-h, x_0+ h\} $
\begin{equation}\label{1D_Six_Order_Scheme}
   u'(x_0) \approx  \alpha_{-2}\;(u_{-2}-u_{2})+\alpha_{-1}\;(u_{-1}-u_{1}) + \alpha_{0}\;u_{0}+\beta_{-1}\;(u'_{-1} +u'_{1}),
\end{equation} 
where $u_i \approx u(x_0+ih)$ and $u'_{i}\approx u'(x_0+ih)$. Then analytical expressions for the weights are obtained by symbolic computation of the linear system
\begin{equation} \label{1D_RBF_HFD_Order6_Matrix}
\scriptsize{
\begin{bmatrix}
1- e^{-16\epsilon^2 h^2} &  e^{-\epsilon^2 h^2}- e^{-9\epsilon^2 h^2} &  e^{-4\epsilon^2 h^2} &  -2h\epsilon^2 e^{-\epsilon^2 h^2}-6h\epsilon^2 e^{-9\epsilon^2 h^2} & 1\\\\
 e^{-\epsilon^2 h^2}- e^{-9\epsilon^2 h^2} & 1-e^{-4\epsilon^2 h^2} & e^{-\epsilon^2 h^2} & -4h\epsilon^2 e^{-4\epsilon^2 h^2} & 1 \\\\
 0 & 0 & 1 & 0 & 1 \\\\
 -2h\epsilon^2 e^{-\epsilon^2 h^2}-6h\epsilon^2 e^{-9\epsilon^2 h^2} & -4h\epsilon^2 e^{-4\epsilon^2 h^2} & 2h\epsilon^2 e^{-\epsilon^2 h^2} & 2\epsilon^2 + e^{-4\epsilon^2 h^2}{(-16\epsilon^4 h^2+2\epsilon^2)} & 0 \\\\
 0 & 0 & 1 & 0 & 0
\end{bmatrix}
\begin{bmatrix}
    \alpha_{-2} \\ \\ \alpha_{-1} \\\\ \alpha_{0} \\\\ \beta_{-1}\\\\  \gamma 
\end{bmatrix}
=
\begin{bmatrix}
    -4h\epsilon^2 e^{-4\epsilon^2 h^2} \\\\
    -2h\epsilon^2 e^{-\epsilon^2 h^2} \\\\
    0 \\\\
    e^{-\epsilon^2 h^2}{(-4\epsilon^4 h^2+2\epsilon^2)} \\\\
   0
\end{bmatrix}
}
\end{equation}
\noindent On substituting the analytical weights in the local truncation error formula 
\begin{eqnarray} \label{1D_Six_Order_Scheme_LTE}
\tau_{0} & = & \; \alpha_{-2} \;(u(x_0-2h)-u(x_0+2h)) + \alpha_{-1}\;( u(x_0-h) - u(x_0+h)) + \alpha_{0}\; u(x_0) \nonumber \\ 
&  & + \beta_{-1} \; (u'(x_0-h) + u'(x_0+h)) - u'(x_0)
\end{eqnarray}  
and making use of Taylor expansions of $u$ and $u'$ values about $x_0$, leading order expression for local truncation error at $x_0$ is obtained in powers of $\epsilon h\; (\epsilon h \ll 1)$. Analytical expressions for formula (\ref{1D_Six_Order_Scheme}) weights and LTE expression are reported in Table \ref{Table1_RBF_HFD}.
\subsection{Eighth order formula}
\noindent Satyanarayana {\it et al.} \cite{Satya_RBF_HFD} proposed RBF-HFD (order 8) formula on a stencil with local support sets  $S=\{ x_0 -2h, x_0-h, x_0, x_0+h, x_0+2h  \} $ and $\mu=\{  x_0-2h, x_0-h, x_0+ h, x_0+2h  \}  $  
\begin{equation}\label{1D_Eighth_Order_Scheme}
       u'(x_0) \approx  \alpha_{-2}\;(u_{-2}-u_{2})+\alpha_{-1}\;(u_{-1}-u_{1}) + \alpha_{0}\;u_{0} +\beta_{-2}\;(u'_{-2}+u'_{2})+\beta_{-1}\;(u'_{-1} +u'_{1}).  
\end{equation}
\noindent The corresponding linear system is symbolically solved for deriving the analytical expressions for the weights in formula (\ref{1D_Eighth_Order_Scheme}).
Then leading order analytical expression for local truncation error ($\tau_0$) for the approximation in (\ref{1D_Eighth_Order_Scheme}) is obtained on substituting the analytical weights and performing Taylor expansions on $u$ and $u'$ about $x_0$.\\
\subsection{Tenth order formula}
\noindent  On a stencil with local support sets 
$S=\{x_0 -3h, x_0 -2h, x_0-h, x_0, x_0+h, x_0+2h, x_0 +3h  \} $ and $\mu=\{  x_0-2h, x_0-h, x_0+ h, x_0+2h  \} $, 
RBF-HFD (order 10) formula \cite{Satya_RBF_HFD} for approximating first derivative is defined
\begin{eqnarray}\label{1D_Tenth_Order_Scheme} 
 u'(x_0) &\approx &\alpha_{-3}(u_{-3}+u_{3})+\alpha_{-2}(u_{-2}+u_{2})+\alpha_{-1}(u_{-1}+u_{1}) +\alpha_{0}\;u_{0}\nonumber \\
 && +\beta_{-2}(u'_{-2}+u'_{2})+\beta_{-1}(u'_{-1} + u'_{1}). 
\end{eqnarray}
The corresponding linear system for the weights in formula (\ref{1D_Tenth_Order_Scheme}) is not directly amenable to symbolic computation in finite time. Therefore, we replace each non-constant term with its truncated Taylor series expansion with 20 terms. Then the analytical expression for the weights are obtained by solving the approximate linear system in {\it Mathematica}.\\~ \\
\noindent The analytical expressions for weights and local truncation errors for RBF-HFD formulas (order 4, 6, 8, and 10) are reported in Table \ref{Table1_RBF_HFD}. Fourth order formula weights and LTEs (Bayona {\it et al.} \cite{bayona2012gaussian}) are included in Table \ref{Table1_RBF_HFD} for the sake of comparison.
It has been observed that the weights of RBF-HFD approximation formulas for first derivative converge to the weights of respective order central compact FD schemes (Collatz \cite{collatznumerical} and Lele \cite{lele1992compact}) in the limit of $\epsilon \rightarrow 0$. Wright and Fornberg have made a similar kind of observation \cite{wright2006scattered} for MQ based RBF-HFD (order 4 and 6), using {\it Mathematica} and Contour-Pad$\acute{e} $ algorithm. The weights of RBF-HFD (order 8 and 10) formulas in the limit of $\epsilon  \rightarrow 0 $ are schematically displayed in (\ref{1D_Eigth_Order_weights_eps0}) and (\ref{1D_Tenth_Order_weights_eps0}), respectively.

\begin{equation}\label{1D_Eigth_Order_weights_eps0}
     \begin{tikzpicture}
     \node at (-3.25,0) [] (c1) {$ u'(x_0) \thickapprox\; \; $};
        \node at (-2,0) [rectangle,draw] (c2) {$\frac{-1}{36}$};
        \node at (0,0) [rectangle, draw] (c3) {$\frac{-4}{9}$};
         \node at (4,0) [rectangle, draw] (c4) {$\frac{-4}{9}$};
          \node at (6,0) [rectangle, draw] (c5) {$\frac{-1}{36}$};
          \node at (7,0) [] (c6) {$u'~~+$};
          \node at (-2,-1) [rectangle,draw] (c7) {$\frac{-25}{216}$};
          \node at (0,-1) [rectangle,draw] (c8) {$\frac{-20}{27}$};
          \node at (2,-1) [rectangle,draw] (c9) {$0$};
        \node at (4,-1) [rectangle,draw] (c10) {$\frac{20}{27}$}; 
        \node at (6,-1) [rectangle,draw] (c11) {$\frac{25}{216}$};
        \node at (6.75,-1) [] (c12) {${\displaystyle\frac{u}{h}}$};
        \draw (c2)--(c3);
        \draw (c3)--(c4);
        \draw (c4)--(c5);
        \draw (c5)--(c6);
        \draw (c7)--(c8);
        \draw (c8)--(c9);
        \draw (c9)--(c10);
        \draw (c10)--(c11);
\end{tikzpicture}
\end{equation}
\begin{equation}\label{1D_Tenth_Order_weights_eps0}
     \begin{tikzpicture}
     \node at (-2,0) [] (c1) {$ u'(x_0) \thickapprox\; $};
        \node at (0,0) [rectangle,draw] (c2) {$\frac{-1}{20}$};
        \node at (2,0) [rectangle,draw] (c3) {$\frac{-1}{2}$};
         \node at (5,0) [rectangle,draw] (c4) {$\frac{-1}{2}$};
          \node at (7,0) [rectangle,draw] (c5) {$\frac{-1}{20}$};
          \node at (8,0) [] (c6) {$u'\;\;+$};
          \node at (-2,-1) [rectangle,draw] (c7) {$\frac{-1}{600}$};
        \node at (0,-1) [rectangle,draw] (c8) {$\frac{-101}{600}$};
         \node at (1.75,-1) [rectangle,draw] (c9) {$\frac{-17}{24}$};
          \node at (3.25,-1) [rectangle,draw] (c10) {$0$};
        \node at (5,-1) [rectangle,draw] (c11) {$\frac{17}{24}$}; 
        \node at (6.75,-1) [rectangle,draw] (c12) {$\frac{101}{600}$};
         \node at (8.5,-1) [rectangle,draw] (c13) {$\frac{1}{600}$};
        \node at (9.25,-1) [] (c14) {${\displaystyle\frac{u}{h}}$};
        \draw (c2)--(c3);
        \draw (c3)--(c4);
        \draw (c4)--(c5);
        \draw (c7)--(c8);
        \draw (c8)--(c9);
        \draw (c9)--(c10);
        \draw (c10)--(c11);
         \draw (c11)--(c12);
          \draw (c12)--(c13);         
\end{tikzpicture}
\end{equation}
\begin{figure}[h!]
    \begin{center}
    \includegraphics[width=100ex,height=50ex]{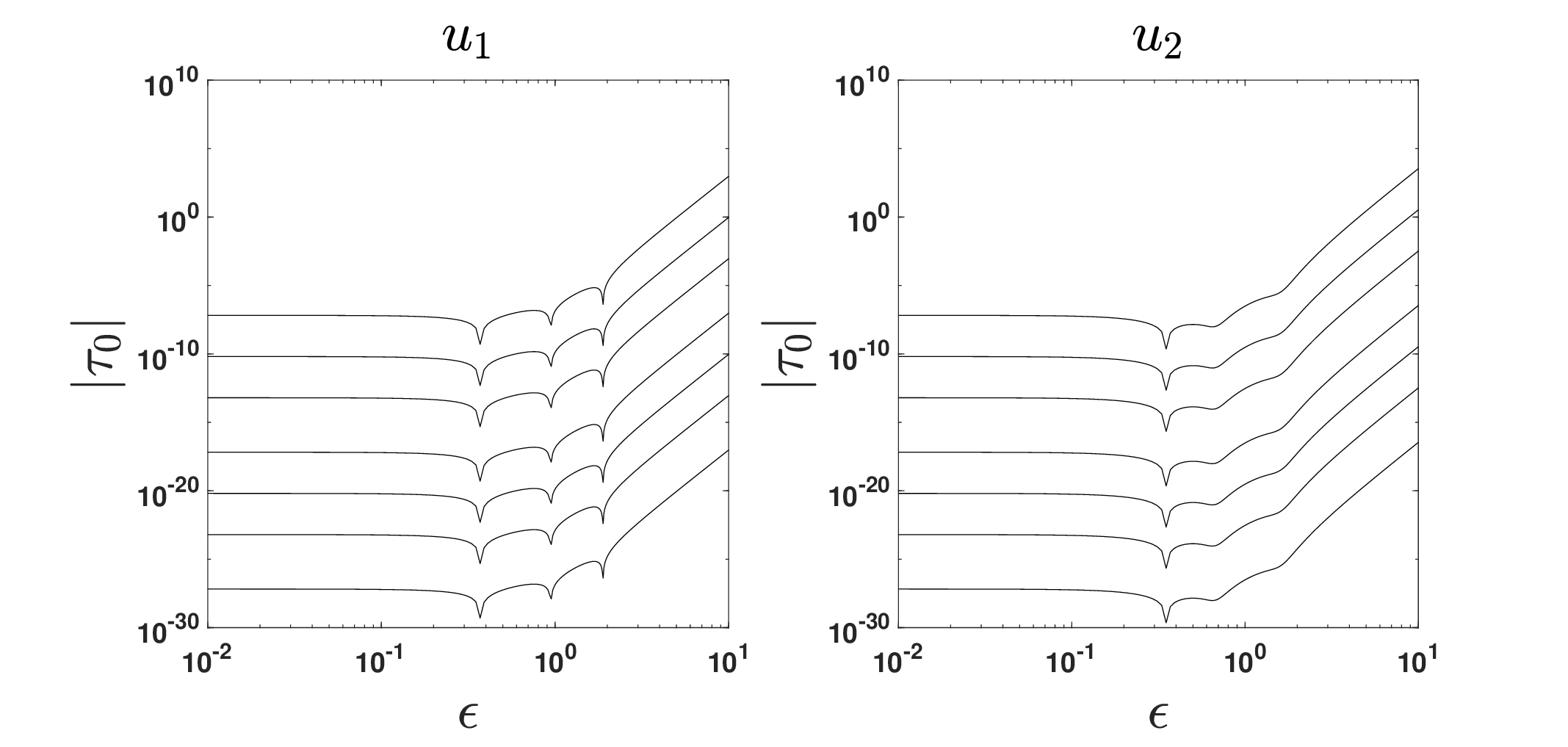}
    \caption{Shape parameter ($\epsilon$) dependence on absolute LTEs ($|\tau_{0}|$) of RBF-HFD (order 10) formula for first derivative approximation of test functions $u_{1}$ and $u_{2}$ for $ h = 0.2,\;0.1,\; 0.05,\; 0.02,\;0.01,\;0.005$\;and\;$ 0.002$, from top to bottom.} 
    \label{fig:1D_Order10_LTE_EPS}
    \end{center}
\end{figure} \\
\noindent In Figure \ref{fig:1D_Order10_LTE_EPS}, we plot the absolute LTEs associated to tenth order RBF-HFD formula approximating first derivative of test functions defined in (\ref{u1}) and (\ref{u2}) against the shape parameter ($\epsilon$). It is observed that the LTE is has a global minimum at some (optimal) value of $\epsilon$. The optimal value is independent of step size ($h$) but depends on the choice of test function and reference point. This trend is seen in Figure \ref{fig:1D_Order10_LTE_EPS} for different step sizes and a range of shape parameter values. Absolute LTEs of sixth and eighth order formulas also display similar trend as the LTE of tenth order formula.\\~ \\
\noindent We now compare the approximation accuracy of GA based RBF-HFD formulas (order 4, 6, 8, 10) with corresponding order MQ based RBF-HFD formulas \cite{Satya_RBF_HFD} and compact FD schemes \cite{collatznumerical, lele1992compact} for first derivative of test functions defined in (\ref{u1})--(\ref{u2}) and fixed value of step size ($h=0.01$). Figure \ref{fig:Comparision_GA_vs_MQ_1D_Order46810_u1_LTE_EPS} depict the variation in absolute LTEs of the formulas against shape parameter ($\epsilon$). In this figure horizontal dashed ($--$) lines represent the absolute LTEs of compact FD (order 4, 6, 8 and 10). The solid (---) curves represent the absolute LTEs of GA based RBF-HFD (order 4, 6, 8 and 10) formulas. The dashed-dot ($-\cdot -$) curves represent the absolute LTEs of respective order MQ based RBF-HFD formulas. It is evident from these figures that in general, the GA based RBF-HFD formulas produce better accuracy as compared to respective order compact FD schemes and MQ based RBF-HFD formulas for a range of shape parameter values. 
\begin{figure}[h!]
    \begin{center}
    \includegraphics[width=45ex,height=45ex]{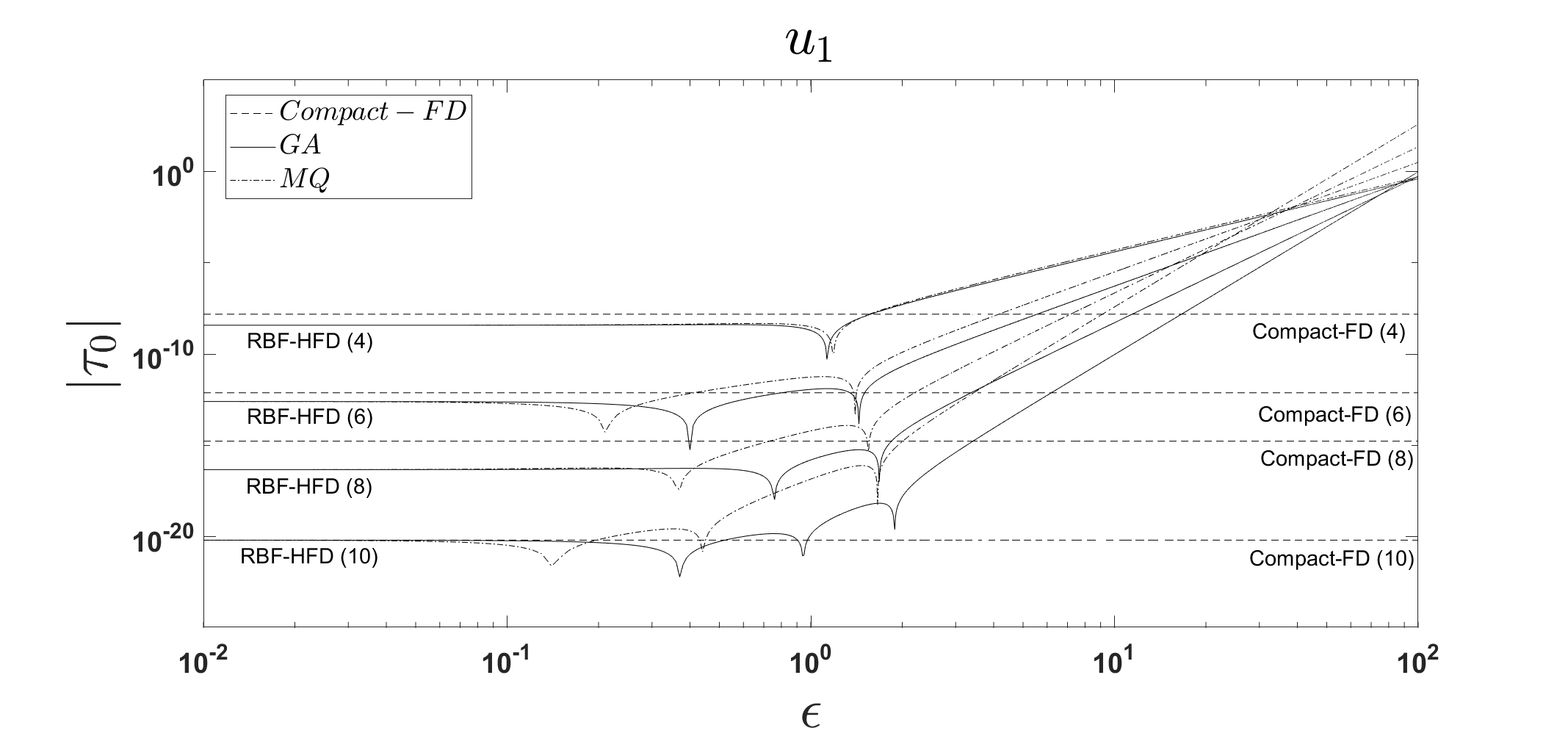} \includegraphics[width=45ex,height=45ex]{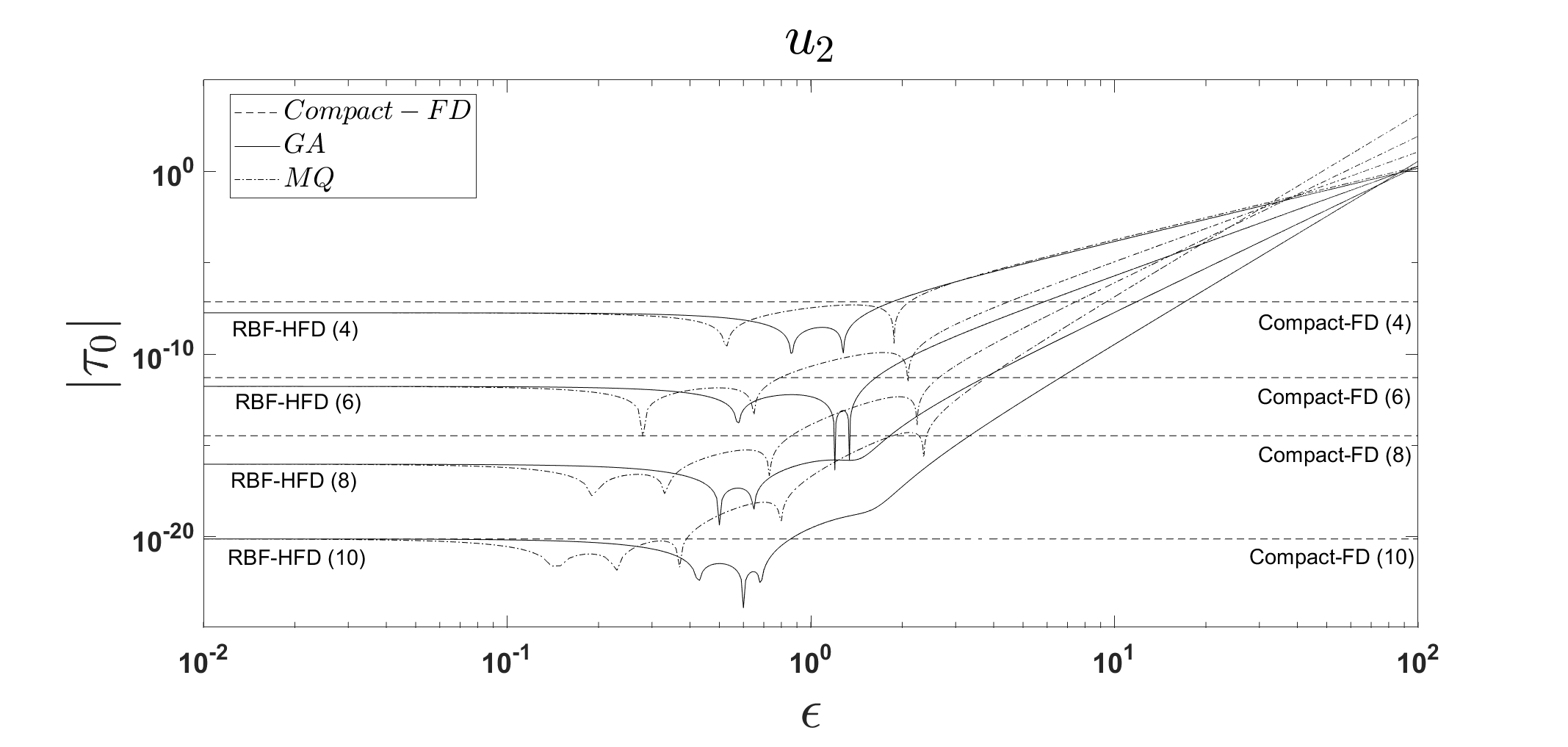}
    \caption{Comparison of shape parameter ($\epsilon$) dependence on absolute LTEs ($|\tau_{0}|$) of MQ \cite{Satya_RBF_HFD} and GA based RBF-HFD (order 4, 6, 8 and 10) and respective order compact FD schemes \cite{collatznumerical, lele1992compact} for first derivative approximation of test function $u_{1}$ at $x_0=0.4$ and $u_{2}$  at $x_0=0.25$.} 
    \label{fig:Comparision_GA_vs_MQ_1D_Order46810_u1_LTE_EPS}
    \end{center}
\end{figure}
 \begin{landscape}
   \begin{table}[t]
    \caption{First derivative approximations using RBF-HFD (order 4,6,8,10) and validation with \cite{collatznumerical}\; (Table 3, P. 538)}
\begin{center} 
\footnotesize{
\begin{tabular}{|c|c|c|} 
 \hline
  & & \\
Stencil & Scheme, weights and local truncation error(LTE)  & Weights and LTE  \\
    & & \\
\hline
& $ u'(x_0)\approx \alpha_{-1}u_{-1} + \alpha_{0}u_{0} + \alpha_{1}u_{1} +\beta_{-1}u'_{-1} + \beta_{1}u'_{1} $   & $\alpha_{-1}=-\alpha_{1}= -\frac{3}{4h},$(\cite{collatznumerical},  T. 3, P. 538),\cite{wright2006scattered}\\
Order 4 \cite{bayona2012gaussian} &$\alpha_{-1}=-\alpha_{1}=\dfrac{2\epsilon^2 h e^{3\epsilon^2 h^2}(-4\epsilon^2 h^2+ e^{4\epsilon^2 h^2}-1)}{-8\epsilon^2 h^2 e^{4\epsilon^2 h^2}+e^{8\epsilon^2 h^2}-1}$ &$ \alpha_{0}=0,$ \\
$x_{-1}\;\;\;\;\;x_{0}\;\;\;\;\;~ ~ x_{1}$&$\alpha_{0}=0,~\beta_{-1}=\beta_{1}=\dfrac{e^{\epsilon^2 h^2}(2\epsilon^2 h^2\cosh({2\epsilon^2 h^2})-\sinh({2\epsilon^2 h^2}))}{4\epsilon^2 h^2-\sinh({4\epsilon^2 h^2})}$&$\beta_{-1}=\beta_{1}=-\frac{1}{4},$\\
\begin{tikzpicture}
  \draw (0,0) circle (0.25);
  \draw (0,0) circle (0.15);
  \draw[line width=1mm] (0.75,0) circle (0.15);
  \draw (1.5,0) circle (0.25);
  \draw (1.5,0) circle (0.15);
\end{tikzpicture}& $ \tau_{0}=\frac{h^4}{120}(-60 \epsilon^{4} u'(x_0)- 20 \epsilon^{2} u^{(3)}(x_0)- u^{(5)}(x_0)) + O(h^6 P_{3}(\epsilon^2))$&$ \tau_{0} \approx \dfrac{h^4}{30}u^{(5)}(x_0)$\\
& &  \\
\hline 
& & \\
&$ u'(x_0) \approx \alpha_{-2}u_{-2}+\alpha_{-1}u_{-1} + \alpha_{0}u_{0} + \alpha_{1}u_{1}+\alpha_{2}u_{2}+\beta_{-1}u'_{-1} + \beta_{1}u'_{1}$&$\alpha_{-2}=-\alpha_{2} = -\frac{1}{36 h}, $ (\cite{collatznumerical},  T. 3, P. 538),\cite{wright2006scattered}\\
Order 6 &$\alpha_{-2}=-\alpha_{2} = \frac{-1}{36 h}-\frac{\epsilon^2 h}{9}-\frac{10 \epsilon^4 h^3}{63}-\frac{8\epsilon^6 h^5}{189}+ O(\epsilon^8 h^7)$& $\alpha_{0}=0,$\\ 
$x_{-2}\;\;\;\;x_{-1}\;\;\;\;\;x_{0}\;\;\;\;\;x_{1}\;\;\;\;\;x_{2}$&$\alpha_{-1}=-\alpha_{1} = \frac{-7}{9 h}-\frac{\epsilon^2 h}{9}+\frac{43 \epsilon^4 h^3}{126}-\frac{43\epsilon^6 h^5}{378}+O(\epsilon^8 h^7)$&$\alpha_{-1}=-\alpha_{1} = \frac{-7}{9 h},$\\
\begin{tikzpicture}
\draw (0,0) circle (0.15);  
\draw (0.75,0) circle (0.25); 
\draw (0.75,0) circle (0.15);
\draw[line width=1mm] (1.5,0) circle (0.15);
\draw (2.25,0) circle (0.25); 
\draw (2.25,0) circle (0.15);
\draw (3.0,0) circle (0.15);
\end{tikzpicture}&$\alpha_{0}=0,~\beta_{-1}=\beta_{1}= -\frac{1}{3}-\frac{\epsilon^2 h^2}{3} +\frac{\epsilon^4 h^4}{42}+\frac{17\epsilon^6 h^6}{126}+ O(\epsilon^8 h^8)$&$\beta_{-1}=\beta_{1}= -\frac{1}{3},$ \\
&$\tau_{0}=\frac{h^6}{1260}\left(840 \epsilon^6 u'(x_0) + 420 \epsilon^4 u^{(3)}(x_0) + 42\epsilon^2  u^{(5)}(x_0) +  u^{(7)}(x_0)\right)  + O(h^8 P_{4}(\epsilon^2))$&$ \tau_{0} \approx -\dfrac{h^6}{420} u^{(7)}(x_0)$ \\
&&\\
\hline
&&\\
& $ u'(x_0) \approx \alpha_{-2}u_{-2}+\alpha_{-1}u_{-1} + \alpha_{0}u_{0} + \alpha_{1}u_{1} + \alpha_{2}u_{2}+\beta_{-2}u'_{-2}+\beta_{-1}u'_{-1} + \beta_{1}u'_{1} +\beta_{2}u'_{2}$&(\cite{collatznumerical},  T. 3, P. 538),\cite{wright2006scattered}\\
&&\\
Order 8 &$\alpha_{-2}=-\alpha_{2}=\frac{-25}{216 h} - \frac{19\epsilon^2 h}{54}-\frac{139\epsilon^4 h^3}{486}+\frac{3532\epsilon^6h^5}{18711}+O(\epsilon^8 h^7)$&$-\alpha_{-2}=-\alpha_{2} = -\frac{25}{216 h}$,\\ 
  $x_{-2}\;\;\;\;\;\;x_{-1}\;\;\;\;\;x_{0}\;\;\;\; x_{1}\;\;\;\;\;\;\;x_{2}$ & $\alpha_{-1}=-\alpha_{1} = \frac{-20}{27h} + \frac{4 \epsilon^2 h}{27} +\frac{106 \epsilon^4 h^3}{243}-\frac{62 \epsilon^6 h^5}{243}+ O(\epsilon^8 h^7)$&$ \alpha_{-1}=-\alpha_{1} =- \frac{20}{27h}$ \\
\begin{tikzpicture}
\draw (0,0) circle (0.15);
\draw (0,0) circle (0.25);
\draw (0.75,0) circle (0.25);
\draw (0.75,0) circle (0.15);
\draw[line width=1mm] (1.5,0) circle (0.15);
\draw (2.25,0) circle (0.15);
\draw (2.25,0) circle (0.25);
\draw (3.0,0) circle (0.15);
\draw (3.0,0) circle (0.25);
\end{tikzpicture}& $\alpha_{0}=0,\beta_{-2}=\beta_{2}= -\frac{1}{36} -\frac{1}{9}\epsilon^2h^2-\frac{13}{81}\epsilon^4h^4-\frac{4}{81}\epsilon^6h^6+O(\epsilon^8 h^8)$ &$\beta_{-1}=\beta_{1}= \frac{-4}{9}$\\
&$\beta_{-1}=\beta_{1}= \frac{-4}{9}-\frac{-4}{9}\epsilon^2h^2+\frac{2}{81}\epsilon^4h^4+\frac{14}{81}\epsilon^6h^6 +O(\epsilon^8 h^8)$&$\alpha_{0}=0,~\beta_{-2}=\beta_{2}=-\frac{1}{36}$\\
&$\tau_{0}= h^8\left(\frac{-2}{3}\epsilon^8 u'(x_0)-\frac{4}{9} \epsilon^6 u^{(3)}(x_0)- \frac{1}{15}\epsilon^4 u^{(5)}(x_0)-\frac{1}{315} \epsilon^2  u^{(7)}(x_0)- \frac{1}{22680}u^{(9)}(x_0)\right)  + O(h^{10} P_{5}(\epsilon^2))$&$\tau_{0} \approx -\dfrac{h^8}{630} u^{(9)}(x_0)$\\
& & \\
\hline
&&\\
& $ u'(x_0) \approx \alpha_{-3}u_{-3}+\alpha_{-2}u_{-2}+\alpha_{-1}u_{-1} + \alpha_{0}u_{0} + \alpha_{1}u_{1} + \alpha_{2}u_{2}+\alpha_{3}u_{3}+\beta_{-2}u'_{-2}+\beta_{-1}u'_{-1} + \beta_{1}u'_{1} +\beta_{2}u'_{2}$&( \cite{lele1992compact} Table 2, P.19)\\
&$\alpha_{-3}=-\alpha_{3}=\frac{-1}{600 h} - \frac{3\epsilon^2 h}{200}+\frac{7639616774597\epsilon^4 h^3}{58982400000}+\frac{6690462911023\epsilon^6h^5}{3276800000}+O(\epsilon^8 h^7)$&$\alpha_{-3}=-\alpha_{3} = -\frac{1}{600 h}$,\\ 
Order 10 &$\alpha_{-2}=-\alpha_{2}=\frac{-101}{600 h} - \frac{71\epsilon^2 h}{150}+\frac{542622915796387\epsilon^4 h^3}{132710400000}+\frac{4920524576353819\epsilon^6h^5}{132710400000}+O(\epsilon^8 h^7)$&$\alpha_{-2}=-\alpha_{2}=\frac{-101}{600 h},$\\
 & $\alpha_{-1}=-\alpha_{1} = \frac{-17}{24h} + \frac{7 \epsilon^2 h}{24} -\frac{53491473198179 \epsilon^4 h^3}{21233664000}-\frac{299239342744897 \epsilon^6 h^5}{10616832000}+ O(\epsilon^8 h^7)$&$\alpha_{-1}=-\alpha_{1} = -\frac{17}{24 h}$,\\
 $x_{-3}\;\;\;\;\;x_{-2}\;\;\;\;\;x_{-1}\;\;\;\;x_{0}\;\;\;\;\;x_{1}\;\;\;\;\;x_{2}\;\;\;\;\;x_{3}$ & $\alpha_{0}=0,~\beta_{-2}=\beta_{2}= -\frac{1}{20}-\frac{\epsilon^2 h^2}{5}+\frac{7641828614597 \epsilon^4 h^4}{4423680000}+\frac{82219926605429\epsilon^6 h^6}{4423680000}+O(\epsilon^8 h^8)$&$\beta_{-1}=\beta_{1}= -\frac{1}{2},\alpha_{0}=0,$ \\
\begin{tikzpicture}
\draw (-0.75,0) circle (0.15);
\draw (0,0) circle (0.15);
\draw (0,0) circle (0.25);
\draw (0.75,0) circle (0.25);
\draw (0.75,0) circle (0.15);
\draw[line width=1mm] (1.5,0) circle (0.15);
\draw (2.25,0) circle (0.15);
\draw (2.25,0) circle (0.25);
\draw (3.0,0) circle (0.15);
\draw (3.0,0) circle (0.25);
\draw (3.75,0) circle (0.15);
\end{tikzpicture}&$\beta_{-1}=\beta_{1}=-\frac{1}{2}- \frac{\epsilon^2 h^2}{2}+\frac{7643155718597 \epsilon^4 h^4}{1769472000}+\frac{29645672092819 \epsilon^6 h^6}{884736000}+O(\epsilon^8 h^8)$&$\beta_{-2}=-\beta_{2} =- \frac{1}{20},$\\
&$\tau_{0} = {h^{10}}(\frac{6}{5}\epsilon^{10} u'(x_0)+\epsilon^8 u^{(3)}(x_0)+\frac{55440}{277200} \epsilon^6 u^{(5)}(x_0)+\frac{3960}{277200} \epsilon^4 u^{(7)}(x_0)+\frac{110}{277200} \epsilon^2 u^{(9)}(x_0)+\frac{1}{277200} u^{(11)}(x_0)))$&$\tau_{0} \approx \frac{144 h^{10}}{39916800} u^{(11)}(x_0)$\\
& $+ O(h^{12} P_{6}(\epsilon^2)$ &\\
 \hline
\end{tabular}
}
\label{Table1_RBF_HFD}
\end{center}
\end{table}
\end{landscape}
\section{Second Derivative Approximations} \label{Second_derivative_approximation}
\noindent Bayona {\it et al.} \cite{bayona2012gaussian} derived analytical expressions of weights and local truncation error for Gaussian based fourth order RBF-HFD formula for second derivative. We derive analytical expressions of weights and LTEs for sixth, eighth and tenth order GA based RBF-HFD formulas for second derivative. The formula stencils are shown in Table \ref{Table2_RBF_HFD}. In these stencils, single circle indicates the function value, whereas, double circle indicates function value and its second derivative value at the node. The obtained formulas converge to respective order compact FD formulas \cite{collatznumerical,lele1992compact} in the limit $\epsilon \rightarrow 0$. 
\begin{figure}[h!]
    \begin{center}
    \includegraphics[width=100ex,height=50ex]{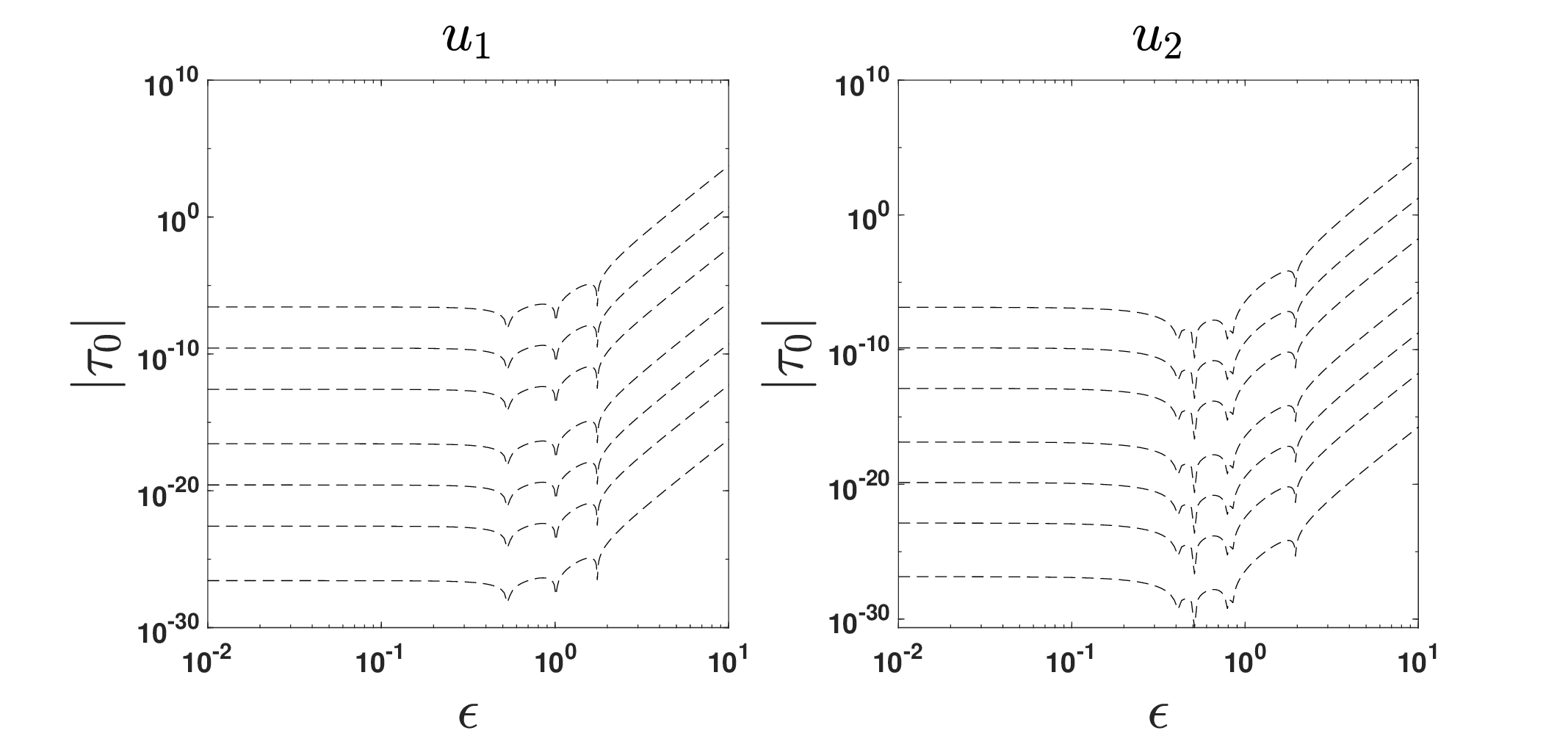}
    \caption{Shape parameter ($\epsilon$)  dependence on absolute LTEs ($|\tau_{0}|$) of RBF-HFD (order 10) formula for second derivative approximation of test functions $u_{1}$ and $u_{2}$ and $ h = 0.2,\;0.1,\; 0.05,\; 0.02,\;0.01,\;0.005$\;and\;$ 0.002$, from top to bottom.} 
    \label{fig:2D_Order10_LTE_EPS}
    \end{center}
\end{figure}
\subsection{Sixth order formula}
\noindent Sixth order formula stencil for second derivative is based on local support sets $S=\{  x_0-2h, x_0-h, x_0, x_0+h, x_0+2h  \} $ and $\mu=\{   x_0-h, x_0+ h  \}$  
\begin{equation}\label{RBF_HFD6_2D}
     u''(x_0) \approx \alpha_{-2}\;(u_{-2}+u_{2})+\alpha_{-1}\;(u_{-1} +u_{1})+ \alpha_{0}\;u_{0} +\beta_{-1}\;(u''_{-1}+u''_{1}).
\end{equation}
\noindent Due to uniform distribution of stencil nodes and central difference approximation of second derivative, number of unknown weights in formula (\ref{RBF_HFD6_2D}) is reduced to four. The weights are obtained by symbolic computation of the corresponding linear system. Expression for local truncation error ($\tau_0$) for the approximation in (\ref{RBF_HFD6_2D}) is obtained on substituting the weights. \\
\subsection{Eighth order formula}
\noindent Eighth order GA based formula \cite{Satya_RBF_HFD} with stencil based on local support sets $S=\{  x_0-2h, x_0-h, x_0, x_0+h, x_0+2h  \} $ and $\mu=\{ x_0  -2h, x_0-h,  x_0+h, x_0+2h  \} $ is proposed as
\begin{equation}\label{RBF_HFD8_2D}
     u''(x_0) \approx \alpha_{-2}(u_{-2}+u_{2})+\alpha_{-1}(u_{-1}+u_{1}) + \alpha_{0}u_{0} +\beta_{-2}(u''_{-2}+u''_{2})+
    \beta_{-1}(u''_{-1} + u''_{1}). 
\end{equation}
\noindent The unknown weights in formula (\ref{RBF_HFD8_2D}) are obtained by solving the corresponding linear system. We have observed that the linear system is directly not amenable to symbolic computation in {\it Mathematica}. Therefore, first we expand each element of the system using Taylor series in powers of $\epsilon h\; (\epsilon h \ll 1 )  $ up to order eighteen and then solve the approximate linear system using {\it Mathematica}.  Expression for local truncation ($\tau_0$) is derived on substituting the analytical weights and performing Taylor expansions about $x_0$. \\
\subsection{Tenth order formula}
\begin{figure}
    \begin{center}
    \includegraphics[width=48ex,height=40ex]{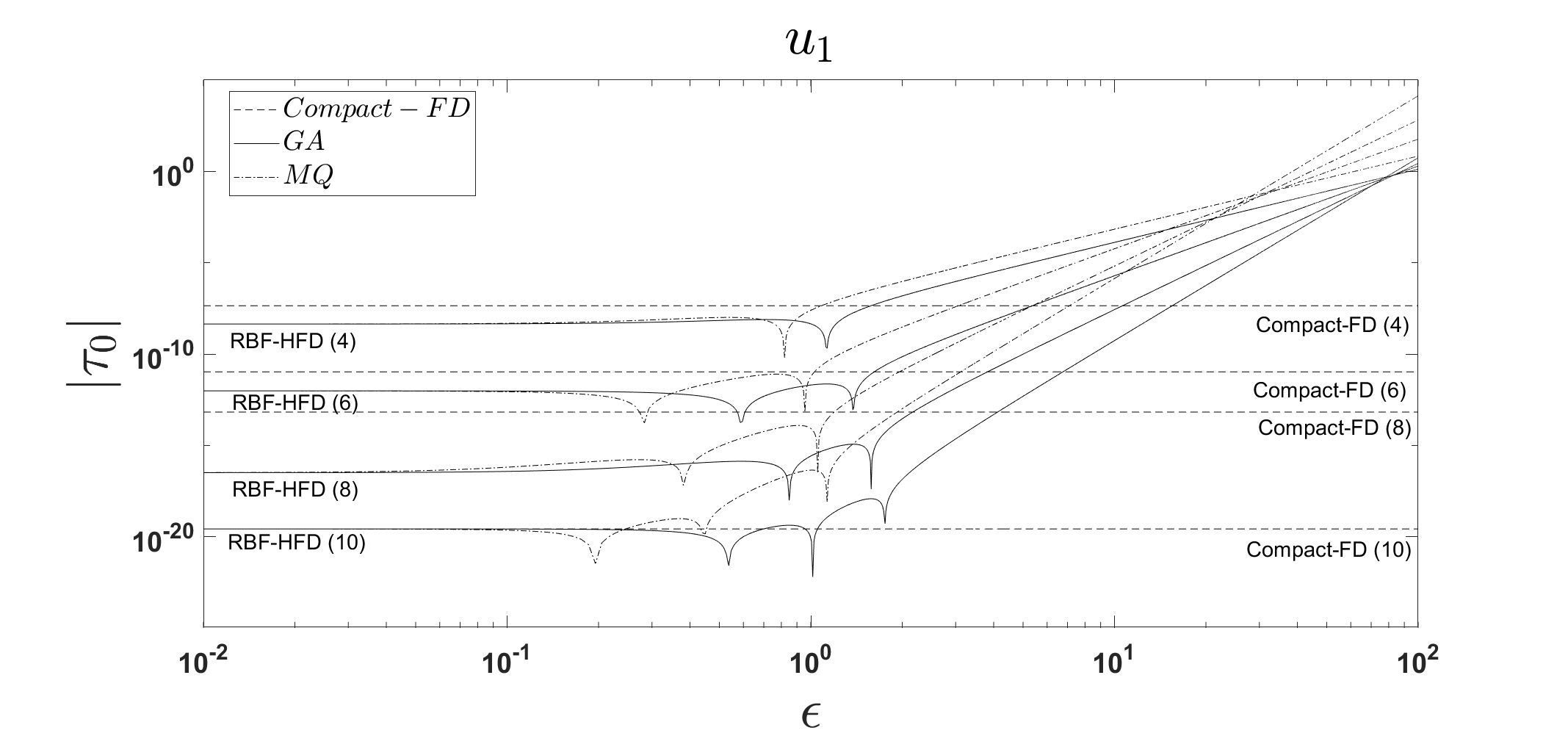}\includegraphics[width=48ex,height=40ex]{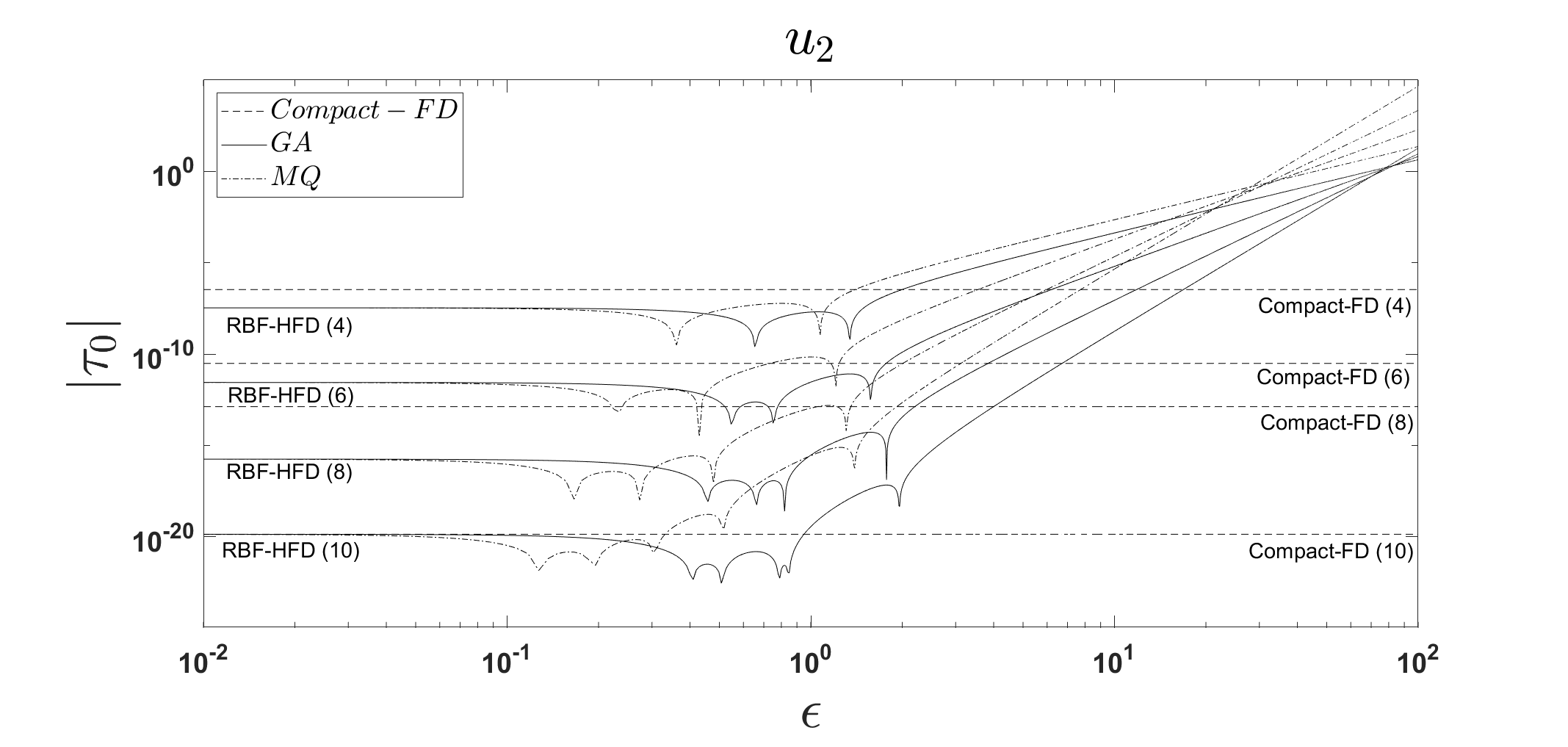}
    \caption{Comparison of shape parameter ($\epsilon$) dependence on absolute LTEs ($|\tau_{0}|$) of MQ \cite{Satya_RBF_HFD} and GA based RBF-HFD (order 4, 6, 8 and 10) and respective order compact FD schemes \cite{collatznumerical, lele1992compact} for second derivative approximations of test function $u_{1}$ at $x_0=0.4$ and $u_{2}$ at $x_0=0.25$.} \label{fig:Comparision_GA_vs_MQ_2D_Order46810_u1_LTE_EPS}
    \end{center}
 \end{figure}

\noindent Approximation of second derivative using tenth order RBF-HFD formula \cite{Satya_RBF_HFD} based on local support sets $S=\{  x_0-3h, x_0-2h, x_0-h, x_0, x_0+h, x_0+2h, x_0+3h   \} $ and $\mu=\{  x_0 -2h, x_0-h, x_0+ h, x_0+2h  \} $ is defined as 
\begin{eqnarray}\label{RBF_HFD10_2D} 
 u''(x_0) &\approx &\alpha_{-3}(u_{-3}+u_{3})+\alpha_{-2}(u_{-2}+u_{2})+\alpha_{-1}(u_{-1}+u_{1}) +\alpha_{0}\;u_{0}\nonumber \\
 && +\beta_{-2}(u''_{-2}+u''_{2})+\beta_{-1}(u''_{-1} + u''_{1}). 
\end{eqnarray}
\noindent The weights in formula (\ref{RBF_HFD10_2D}) are obtained by solving the corresponding linear system. Again this linear system is also directly not solvable in {\it Mathematica}. Therefore, first we expand each element of the system using Taylor series in powers of $\epsilon h\; (\epsilon h \ll 1 ) $ up to order twenty four and then solve the approximate linear system using {\it Mathematica}. Expression for local truncation error ($\tau_0$) for the approximation in (\ref{RBF_HFD10_2D}) is derived on substituting the weights and making use of Taylor series expansion of $u$ and $u''$ values about the reference point.\\~ \\
\noindent The analytical expressions of weights and local truncation errors for sixth order formula (\ref{RBF_HFD6_2D}),  eighth order formula (\ref{RBF_HFD8_2D}) and tenth order formula (\ref{RBF_HFD10_2D}) are reported in the second, third and fourth rows of Table \ref{Table2_RBF_HFD}, respectively. For the sake of comparison, expressions for fourth order formula weights and local truncation error are reported in the first row of Table  \ref{Table2_RBF_HFD}. The obtained formula weights in the limit of $\epsilon$ tending to zero converge to respective order compact FD formula weights  (order 4, 6, 8) from Collatz \cite{collatznumerical} (Table 3, Page. 538-539) and (order 10) from Lele \cite{lele1992compact} (Table 3, Page. 21). We have included the compact FD weights and LTEs in the last column of Table \ref{Table2_RBF_HFD}. Bayona {\it et al.} \cite{bayona2012gaussian} observed convergence for fourth order Gaussian based RBF-HFD weights in the limit of $\epsilon$ tending to zero. The weights for proposed RBF-HFD (order 8 and 10) formulas, in the limit $\epsilon \rightarrow 0$, are schematically shown in stencils (\ref{2D_RBF_HFD8_Weights_eps0}) and (\ref{2D_RBF_HFD10_Weights_eps0}), respectively. \\
\begin{equation}
     \begin{tikzpicture}\label{2D_RBF_HFD8_Weights_eps0}
     \node at (-2,0) [] (c1) {$ u''(x_0) \thickapprox\; $};
        \node at (0,0) [rectangle,draw] (c2) {$\frac{-23}{2358}$};
        \node at (1.75,0) [rectangle,draw] (c3) {$\frac{-344}{1179}$};
        \node at (6,0) [rectangle,draw] (c4) {$\frac{-344}{1179}$};
          \node at (7.75,0) [rectangle,draw] (c5) {$\frac{-23}{2358}$};
          \node at (9,0) [] (c6) {$u''\;\;+$};
          \node at (0,-1) [rectangle,draw] (c7) {$\frac{155}{786}$};
          \node at (2,-1) [rectangle,draw] (c8) {$\frac{320}{393}$};
          \node at (4,-1) [rectangle,draw] (c9) {$\frac{-265}{131}$};
        \node at (6,-1) [rectangle,draw] (c10) {$\frac{320}{393}$}; 
        \node at (8,-1) [rectangle,draw] (c11) {$\frac{155}{786}$};
        \node at (9,-1) [] (c12) {${\displaystyle\frac{u}{h^2}}$};
        \draw (c2)--(c3);
        \draw (c3)--(c4);
        \draw (c4)--(c5);
        \draw (c7)--(c8);
        \draw (c8)--(c9);
        \draw (c9)--(c10);
        \draw (c10)--(c11);
\end{tikzpicture}
\end{equation}
\begin{equation}
     \begin{tikzpicture}\label{2D_RBF_HFD10_Weights_eps0}
     \node at (-3,0) [] (c1) {$ u''(x_0) \thickapprox$};
        \node at (-1,0) [rectangle, draw] (c2) {$\frac{-43}{1798}$};
        \node at (1,0) [rectangle, draw] (c3) {$\frac{-334}{899}$};
         \node at (5,0) [rectangle, draw] (c4) {$\frac{-334}{899}$};
          \node at (6.75,0) [rectangle, draw] (c5) {$\frac{-43}{1798}$};
          \node at (8,0) [] (c6) {$u''\;\;+$};
    \node at (-3,-1) [rectangle,draw] (c7) {$\frac{-79}{16182}$};
    \node at (-1,-1) [rectangle,draw] (c8) {$\frac{519}{1798}$};
    \node at (1,-1) [rectangle,draw] (c9) {$\frac{1065}{1798}$};
    \node at (3,-1) [rectangle,draw] (c10) {$-\frac{14335}{8091}$};
    \node at (5,-1) [rectangle,draw] (c11) {$\frac{1065}{1798}$}; 
    \node at (7,-1) [rectangle,draw] (c12) {$\frac{519}{1798}$};
    \node at (9,-1) [rectangle,draw] (c13) {$\frac{-79}{16182}$};
        \node at (10,-1) [] (c14) {${\displaystyle\frac{u}{h^2}}$};
        \draw (c2)--(c3);
        \draw (c3)--(c4);
        \draw (c4)--(c5);
        \draw (c7)--(c8);
        \draw (c8)--(c9);
        \draw (c9)--(c10);
        \draw (c10)--(c11);
         \draw (c11)--(c12);
          \draw (c12)--(c13);
\end{tikzpicture}
\end{equation}
\noindent Figure \ref{fig:2D_Order10_LTE_EPS} depict the dependence of shape parameter on local truncation error of GA based RBF-HFD (order 10) formula for second derivative. The absolute LTEs are computed for test functions $u_1$ and $u_2$ against a wide range of shape parameter $\epsilon$ values and various step-sizes. It is clear from these figures that, the shape parameter is independent of step size and depends on the choice of test function and its second derivative values at the reference point.  We now compare the approximation accuracy of Gaussian based RBF-HFD formulas (order 4, 6, 8 and 10) with respective order compact FD schemes and MQ based RBF-HFD formulas \cite{Satya_RBF_HFD}. For this purpose we have chosen two test functions $u_1$, $u_2$ defined in equations (\ref{u1})--(\ref{u2}). Figure \ref{fig:Comparision_GA_vs_MQ_2D_Order46810_u1_LTE_EPS},
correspond to approximation accuracy of the test functions $u_1$, $u_2$ of second derivatives. Where the horizontal dashed lines ($--$) represent the absolute LTEs of compact FD (order 4, 6, 8 and 10) schemes. The solid line curves represent absolute LTEs of Gaussian based RBF-HFD formulas (order 4, 6, 8 and 10), whereas dashed-dot ($-\cdot-$) line curves represent absolute LTEs of MQ based RBF-HFD formulas \cite{Satya_RBF_HFD} (order 4, 6, 8 and 10). It is observed from these figures that, Gaussian based RBF-HFD formulas produces the most accurate results as compared with compact FD formulas and MQ based RBF-HFD formulas for certain values of shape parameter.  
\begin{landscape}
 \begin{table}[t]
    \caption{Second derivative approximations using RBF-HFD (order 4, 6, 8,10) formulas}
\begin{center} 
\footnotesize{
\begin{tabular}{|c|c|c|} 
 \hline
 & & \\
   Stencil  & Scheme, weights and local truncation error (LTE) & Weights and LTE  \\
\hline
& & \\
& $ u''(x_0) = \alpha_{-1}u_{-1} + \alpha_{0}u_{0} + \alpha_{1}u_{1} +\beta_{-1}u''_{-1} + \beta_{1}u''_{1} $ & ( \cite{collatznumerical} Table 3, P. 538-539),\cite{wright2006scattered}\\
&  $\alpha_{-1}=\alpha_{1}= \frac{6}{5h^2} + \frac{42 \epsilon^2}{125} -\frac{3021\epsilon^4 h^2}{6250}- \frac{92969\epsilon^6 h^4}{468750}+O(\epsilon^8 h^6)$ & $\alpha_{-1}=\alpha_{1}=\frac{6}{5h^2}, $ \\ 

Order 4 \cite{bayona2012gaussian} & $ \alpha_{0}=\frac{-12}{5 h^2} -\frac{84 \epsilon^2}{125}+\frac{3021\epsilon^4 h^2}{3125}+\frac{92969 \epsilon^6 h^4}{234375}+O(\epsilon^8 h^6)$ & $\alpha_{0}=-\frac{12}{5 h^2},$\\
$x_{-1}\;\;\; x_{0}\;\;\;\;\;x_{1}$& $ \beta_{-1}=\beta_{1}=\frac{-1}{10}-\frac{26250\epsilon^2 h^2}{156250} -\frac{16925\epsilon^4 h^4}{156250}-\frac{651\epsilon^6 h^6}{156250}+O(\epsilon^8 h^8)$&$\beta_{-1}=\beta_{1}=-\frac{1}{10},$\\
\begin{tikzpicture}
\draw (0,0) circle (0.25);
  \draw (0,0) circle (0.15);
  \draw[line width=1mm] (0.75,0) circle (0.15);
  \draw (1.5,0) circle (0.25);
  \draw (1.5,0) circle (0.15);
\end{tikzpicture}& $ \tau_{0}=\frac{h^4}{200}(-140 \epsilon^{4} u''(x_0)- 28 \epsilon^{2} u^{(4)}(x_0) - u^{(6)}(x_0)) + O(h^6 P_{3}(\epsilon^2))$&  $ \tau_{0} \approx \dfrac{h^4}{90} u^{(6)}(x_i)$\\
& & \\
\hline 
& & \\
& $ u''(x_0) \approx \alpha_{-2}u_{-2}+\alpha_{-1}u_{-1} + \alpha_{0}u_{0} + \alpha_{1}u_{1} + \alpha_{2}u_{2}+\beta_{-1}u''_{-1} + \beta_{1}u''_{1} $& ( \cite{collatznumerical} Table 3, P. 538-539),\cite{wright2006scattered} \\
Order 6  &  $\alpha_{-2}= \alpha_{2} = \frac{3}{44 h^2} + \frac{1035 \epsilon^2}{3388} +\frac{1463277\epsilon^4 h^2}{2608760} + \frac{12613961ep^6 h^4}{28696360}+ O(\epsilon^8 h^6)$  &$\alpha_{-2}= \alpha_{2} = \frac{3}{44 h^2},$ \\
 $x_{-2}\;\; \;\;\;\;x_{-1}\;\;\;\; x_{0}\;\;\;\;\;\; x_{1}\;\;\;\; \; \; x_{2}$& $\alpha_{-1}= \alpha_{1} =\frac{12}{11 h^2} - \frac{414 \epsilon^2}{847}-\frac{538116\epsilon^4 h^2}{326095} -\frac{4542297 \epsilon^6 h^4}{7174090}+ O(\epsilon^8 h^6) $ & $\alpha_{-1}= \alpha_{1} =\frac{12}{11 h^2}$ \\
\begin{tikzpicture}
\draw (0,0) circle (0.15);  
\draw (0.75,0) circle (0.25); 
\draw (0.75,0) circle (0.15);
\draw[line width=1mm] (1.5,0) circle (0.15);
\draw (2.25,0) circle (0.25); 
\draw (2.25,0) circle (0.15);
\draw (3.0,0) circle (0.15);
\end{tikzpicture}& $\alpha_{0}=\frac{-51}{22 h^2}+\frac{621\epsilon^2}{1694} +\frac{2841651 \epsilon^4 h^2}{1304380}+\frac{5555227\epsilon^6 h^4}{14348180}+ O(\epsilon^8 h^6)$& $\alpha_{0}=-\frac{51}{22 h^2},$\\
&$\beta_{-1}=\beta_{1}= -\frac{2}{11}-\frac{79381958040 \epsilon^2 h^2}{216542732560} -\frac{ 64254114540 \epsilon^4 h^4}{216542732560}-\frac{8624761068 \epsilon^6 h^6}{216542732560}+
O(\epsilon^8 h^8)$ & $\beta_{-1}=\beta_{1}= -\frac{2}{11},$\\
& $ \tau_{0} = \frac{-23 h^6}{2520}\left(2520 \epsilon^6 u''(x_i) +756 \epsilon^4 u^{(4)}(x_i) + 54 \epsilon^2  u^{(6)}(x_i) +  u^{(8)}(x_i)  \right)  + O(h^8 P_{4}(\epsilon^2))$& $ \tau_{0} \approx \dfrac{-23 h^6}{5040} u^{(8)}(x_i)$\\
& &\\
\hline
& &\\
& $ u''(x_0) = \alpha_{-2}u_{-2}+\alpha_{-1}u_{-1} + \alpha_{0}u_{0} + \alpha_{1}u_{1} + \alpha_{2}u_{2}+\beta_{-2}u''_{-2}+
    \beta_{-1}u''_{-1} + \beta_{1}u''_{1} +\beta_{2}u''_{2}  $&( \cite{collatznumerical} Table 3, P. 538-539),\cite{wright2006scattered} \\
Order 8 &  $\alpha_{-2}= \alpha_{2} = \frac{155}{786 h^2} - \frac{337172\epsilon^2}{463347}+\frac{981384524105795 \epsilon^4 h^2}{978944714496} + O(\epsilon^6 h^4)$ &$\alpha_{-2}= \alpha_{2} = \frac{155}{786 h^2}$\\ 
$x_{-2}\;\;\; \; \; \; x_{-1}\;\;   x_{0}\;\;\;\; \; x_{1}\;\;\;\; \; \; x_{2}$ & $\alpha_{-1}= \alpha_{1} = \frac{320}{393 h^2} - \frac{723008\epsilon^2}{463347}-\frac{131554452784583 \epsilon^4 h^2}{61184044656} + O(\epsilon^6 h^4) $ &$\alpha_{-1}= \alpha_{1} = \frac{320}{393 h^2}$\\
\begin{tikzpicture}
\draw (0,0) circle (0.15);
\draw (0,0) circle (0.25);
\draw (0.75,0) circle (0.25);
\draw (0.75,0) circle (0.15);
\draw[line width=1mm] (1.5,0) circle (0.15);
\draw (2.25,0) circle (0.15);
\draw (2.25,0) circle (0.25);
\draw (3.0,0) circle (0.15);
\draw (3.0,0) circle (0.25);
\end{tikzpicture}&$\alpha_{0}=\frac{-265}{131 h^2}+\frac{257224\epsilon^2}{154449}+\frac{374495573482511 \epsilon^4 h^2}{163157452416}
+O(\epsilon^6 h^4)$& $\alpha_{0}=\frac{-265}{131 h^2},$\\   
&$\beta_{-2}=\beta_{2}=\frac{-23}{2358}-\frac{76472 \epsilon^2 h^2}{1390041}-\frac{111386667913345\epsilon^4 h^4}{1468417071744}+O(\epsilon^6 h^6)$&$\beta_{-2}=\beta_{2}=\frac{-23}{2358}$\\
&$\beta_{-1}=\beta_{1}=\frac{-344}{1179}-\frac{862048\epsilon^2 h^2}{1390041}-\frac{313528367747261\epsilon^4 h^4}{367104267936}+O(\epsilon^6 h^6)$&$\beta_{-1}=\beta_{1}=\frac{-344}{1179}$\\
&$\tau_{0}=\frac{-79 h^8}{2971080} (55440 \epsilon^8 u''(x_0)+22176 \epsilon^6 u^{(4)}(x_0)+2376 \epsilon^4 u^{(6)}(x_0)+88\epsilon^2  u^{(8)}(x_0)+u^{(10)}(x_0))+ O(h^{10} P_{5}(\epsilon^2))$&$\tau_{0} \approx \frac{79 h^8}{1260} u^{(9)}(x_0)$\\
& &\\
\hline
& $ u''(x_0)=\alpha_{-3}u_{-3}+\alpha_{-2}u_{-2}+\alpha_{-1}u_{-1} + \alpha_{0}u_{0} + \alpha_{1}u_{1} + \alpha_{2}u_{2}+\alpha_{3}u_{3}+\beta_{-2}u''_{-2}+\beta_{-1}u''_{-1} + \beta_{1}u''_{1} +\beta_{2}u''_{2}$&( \cite{lele1992compact} Table 3, P.21)\\
&$\alpha_{-3}= \alpha_{3}=\frac{79}{16182 h^2} +\frac{5270785\epsilon^2 }{106682532}+\frac{4214433321515\epsilon^4 h^2}{18989704061064}+O(\epsilon^6 h^4)$&$\alpha_{-3}= \alpha_{3} = \frac{79}{16182 h^2}$,\\
Order 10 &$\alpha_{-2}= \alpha_{2}=\frac{519}{1798 h^2} + \frac{16456115\epsilon^2}{17780422}+\frac{501398734735\epsilon^4 h^2}{1054983558948}+O(\epsilon^6 h^4)$&$\alpha_{-2}=-\alpha_{2}=\frac{519}{1798 h^2},$\\
& $\alpha_{-1}= \alpha_{1} = \frac{1065}{1798 h^2} - \frac{79866475\epsilon^2}{35560844}-\frac{4626851267795 \epsilon^4 h^2}{2109967117896}+ O(\epsilon^6 h^4)$&$\alpha_{-1}=\alpha_{1} = \frac{1065}{1798 h^2}$,\\
 $x_{-3}\;\; \;\;\;\;x_{-2}\;\;\;\;\; x_{-1}\;\;\;\;x_{0}\;\;\;\; \; x_{1}\;\;\;\;\;x_{2}\;\;\;\;\; \;x_{3}$&$\alpha_{0}=-\frac{14335}{8091 h^2}+\frac{67795975 \epsilon^2}{26670633}+\frac{14201025431705 \epsilon^4 h^2}{4747426015266}+O(\epsilon^6 h^4),$&$\alpha_{0}=-\frac{14335}{8091 h^2},$ \\
\begin{tikzpicture}
\draw (-0.75,0) circle (0.15);
\draw (0,0) circle (0.15);
\draw (0,0) circle (0.25);
\draw (0.75,0) circle (0.25);
\draw (0.75,0) circle (0.15);
\draw[line width=1mm] (1.5,0) circle (0.15);
\draw (2.25,0) circle (0.15);
\draw (2.25,0) circle (0.25);
\draw (3.0,0) circle (0.15);
\draw (3.0,0) circle (0.25);
\draw (3.75,0) circle (0.15);
\end{tikzpicture}&$\beta_{-2}=\beta_{2}= -\frac{43}{1798}-\frac{3822325\epsilon^2 h^2}{26670633}-\frac{578232844745\epsilon^4 h^4}{1582475338422}+O(\epsilon^6 h^6)$&$\beta_{-1}=\beta_{1}= -\frac{334}{899},$\\
&$\beta_{-1}=\beta_{1}=-\frac{334}{899}-\frac{21525725\epsilon^2 h^2}{26670633}-\frac{771306629605\epsilon^4 h^4}{1582475338422}
+O(\epsilon^6 h^6)$&$\beta_{-2}=\beta_{2} = -\frac{43}{1798},$\\
&$\tau_{0}= \frac{619 h^{10}}{29904360}(1441440\epsilon^{10} u''(x_0)+720720 \epsilon^8 u^{(4)}(x_0)+102960 \epsilon^6 u^{(6)}(x_0)+5720 \epsilon^4 u^{(8)}(x_0))$&$\tau_{0} \approx \frac{619 h^{10}}{299043360} u^{(12)}(x_0)$\\
&$+130\epsilon^2 u^{(10)}(x_0)+u^{(12)}(x_0))+O(h^{12} P_{6}(\epsilon^2))$&\\
\hline
\end{tabular}
}
\label{Table2_RBF_HFD}
\end{center}
\end{table}
\end{landscape}
\section{2D-Laplacian Operator Approximations}\label{Two_dimensional_Laplacian_approximation}

\noindent Satyanarayana {\it et al.} \cite{Satya_RBF_HFD} obtained MQ based RBF-HFD (order 4 and 6) formulas to approximate 2D-Laplacian. We derive analytical expressions of weights and LTEs for GA based RBF-HFD (order 4 and 6) formulas. In the limit of $\epsilon$ tending to zero,  the formulas converge to fourth and sixth order compact FD formulas reported in Collatz (\cite{collatznumerical}, p. 542-543). In order to validate the formulas, we consider the following test functions from \cite{wright2006scattered}
\begin{eqnarray}
    u_{4}(x,y)&=&e^{-(x-0.25)^{2}-(y-0.5)^2} \sin(\pi x) \cos(2 \pi y),\; (x_0, y_0) =(0.25, 0.25), \label{Test_function4_Laplacian}  \\
    u_{5}(x,y)&=&\dfrac{25}{25+(x-0.2)^{2}+2 y^2},\; (x_0, y_0) =(0, 0), \label{Test_function5_Laplacian} \\
    u_{6}(x,y)&=&e^{x} \tanh(\dfrac{y}{\sqrt{2}}),\; (x_0, y_0) =(0.1, 0.2), \label{Test_function6_Laplacian}
\end{eqnarray}
and test functions from \cite{ding2005error}
\begin{eqnarray}
u_{7}(x,y)&=& \frac{3}{4}\; e^{-\left(\frac{(9x-2)^2 + (9y-2)^2}{4}\right)}+\frac{3}{4} \; e^{\left(-\frac{(9x+1)^2}{49} -\frac{9y+1}{10}  \right)}+\frac{1}{2} \; e^{-\left(\frac{(9x-7)^2 + (9y-3)^2}{4}\right)} \nonumber \\
 & & -\frac{1}{5} \; e^{-\left((9x-4)^{2}+(9y-7)^{2},\right)},\; (x_0, y_0) =(0.1, 0.2), \label{Test_function7_Laplacian} \\
u_{8}(x,y)&=&\left(1-\frac{x}{2}\right)^{6}\left(1-\frac{y}{2}\right)^{6}+1000(1-x)^3 x^3 (1-y)^3 y^3 + y^6  \left(1-\frac{x}{2}\right)^{6} \nonumber \\
& & +x^6 \left(1-\frac{y}{2}\right)^{6},\;\; (x_0, y_0) =(0.1, 0.2),  \label{Test_function8_Laplacian} \\
u_{9}(x,y)&=&\sin(\pi x) \sin(\pi y),\; (x_0, y_0) =(0.1, 0.2). \label{Test_function9_Laplacian}
\end{eqnarray}

\subsection{Fourth order formula}
\noindent Fourth order RBF-HFD formula for approximating 2D-Laplacian is based on a stencil with local support sets
\begin{eqnarray}
S&=&\{  (x_0, y_0), (x_0+ h,y_0), (x_0,y_0+h), (x_0-h,y_0), (x_0,y_0-h), (x_0+h, y_0+h), \nonumber \\
&& (x_0-h,y_0+ h), (x_0-h, y_0-h), (x_0+h, x_0-h) \}, \nonumber \\
\mu &=&\{(x_0+h,y_0), (x_0,y_0+h), (x_0-h,y_0), (x_0,y_0-h)  \}. \nonumber 
\end{eqnarray}
In view of uniform distribution of stencil nodes and central difference approximation of 2D-Laplacian, RBF-HFD (order 4) formula may be written as
\begin{eqnarray}\label{RBF_HFD4_Scheme_Laplacian}
    \Delta u(x_0,y_0) &\approx &\alpha_{0}\; u_{0,0} + \alpha_{1} ( u_{1,0} +  u_{0,1} +  u_{-1,0} +  u_{0,-1} ) + \alpha_{2} ( u_{1,1} +  u_{-1,1} +  u_{-1,-1} +  u_{1,-1} ) \nonumber \\ 
   & & + \beta_{1} \left( (\Delta u)_{1,0} +  (\Delta u)_{0,1}  +  (\Delta u)_{-1,0}  + (\Delta u)_{0,-1} \right). 
\end{eqnarray}
\noindent Here $u_{i,j}\approx u(x_0+ih,y_0+jh)$ and $(\Delta u)_{i,j}\approx \Delta u(x_0+ih,y_0+jh)$. The unknown weights in formula (\ref{RBF_HFD4_Scheme_Laplacian}) may be obtained by symbolically solving the corresponding linear system. Several elements of this linear system are large expressions involving Gaussian radial basis functions. Therefore, this system is directly not amenable to symbolic computation in {\it Mathematica}. To circumvent this problem, we express each non-constant term of this system using Taylor series expansion in powers of ($\epsilon h $) up to order sixteen. Then we solve the approximate linear system via symbolic computation. On substituting the weights so obtained and making use of Taylor expansions of $u$ and $\Delta u$ terms about $(x_0, y_0)$, we extract leading order expression for local truncation error ($\tau_0$) for the approximation in (\ref{RBF_HFD4_Scheme_Laplacian}). The expressions for weights and local truncation error are reported in the second row of Table \ref{Table3_RBF_HFD}. \\  

\noindent To validate the obtained weights with benchmark results, we have included the weights and local truncation error of classical compact FD (order 4) \cite{collatznumerical} in the last column of Table  \ref{Table3_RBF_HFD}. It is observed that the GA based RBF-HFD (order 4) formula converge to compact FD (order 4) scheme as $ \epsilon \rightarrow 0 $. Wright and Fornberg \cite{wright2006scattered} and Satyanarayana {\it et al.} \cite{Satya_RBF_HFD} also observed convergence for the MQ based formulas in the flat limit. \\
\begin{figure}[h!]
    \begin{center}
     \includegraphics[width=48ex,height=40ex]{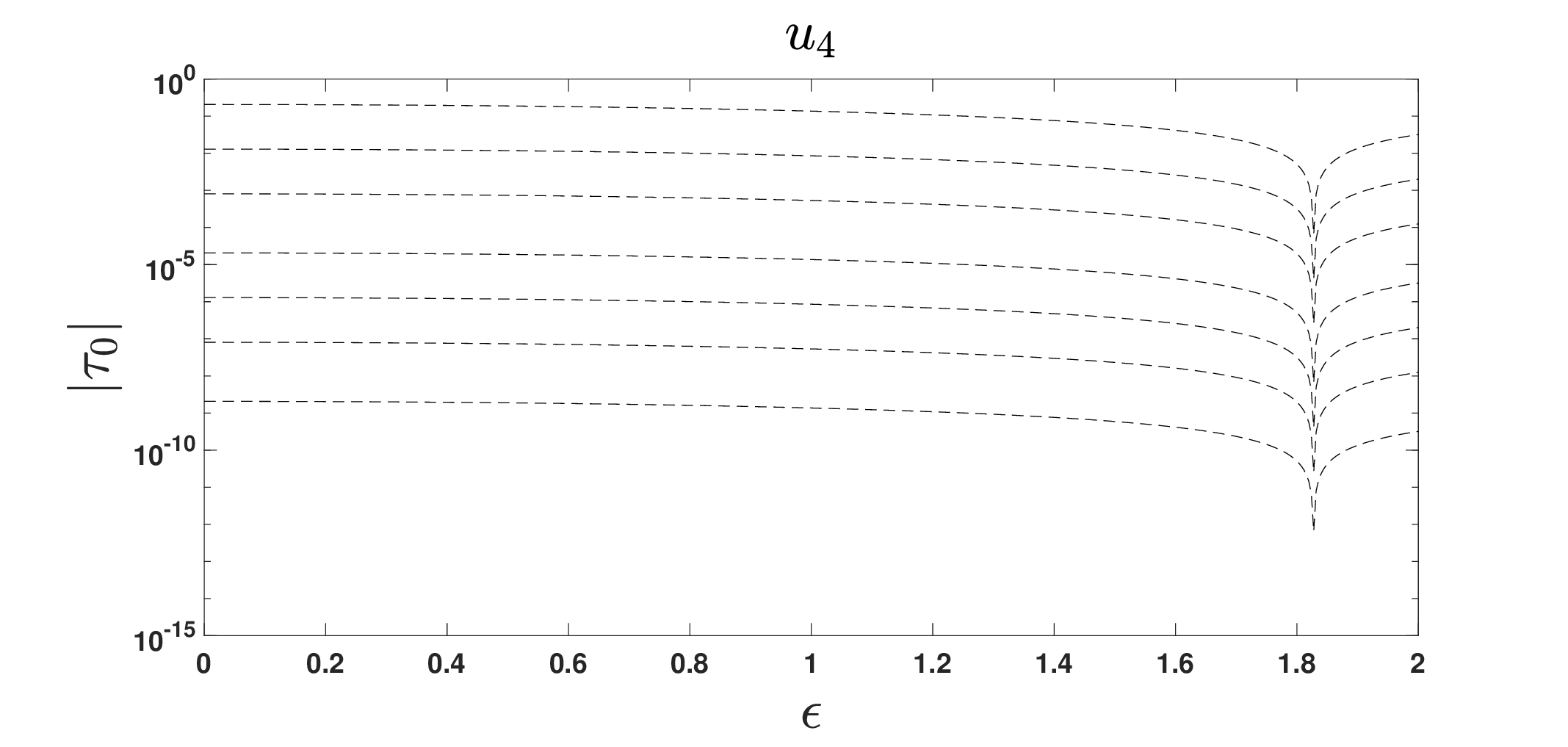}\includegraphics[width=48ex,height=40ex]{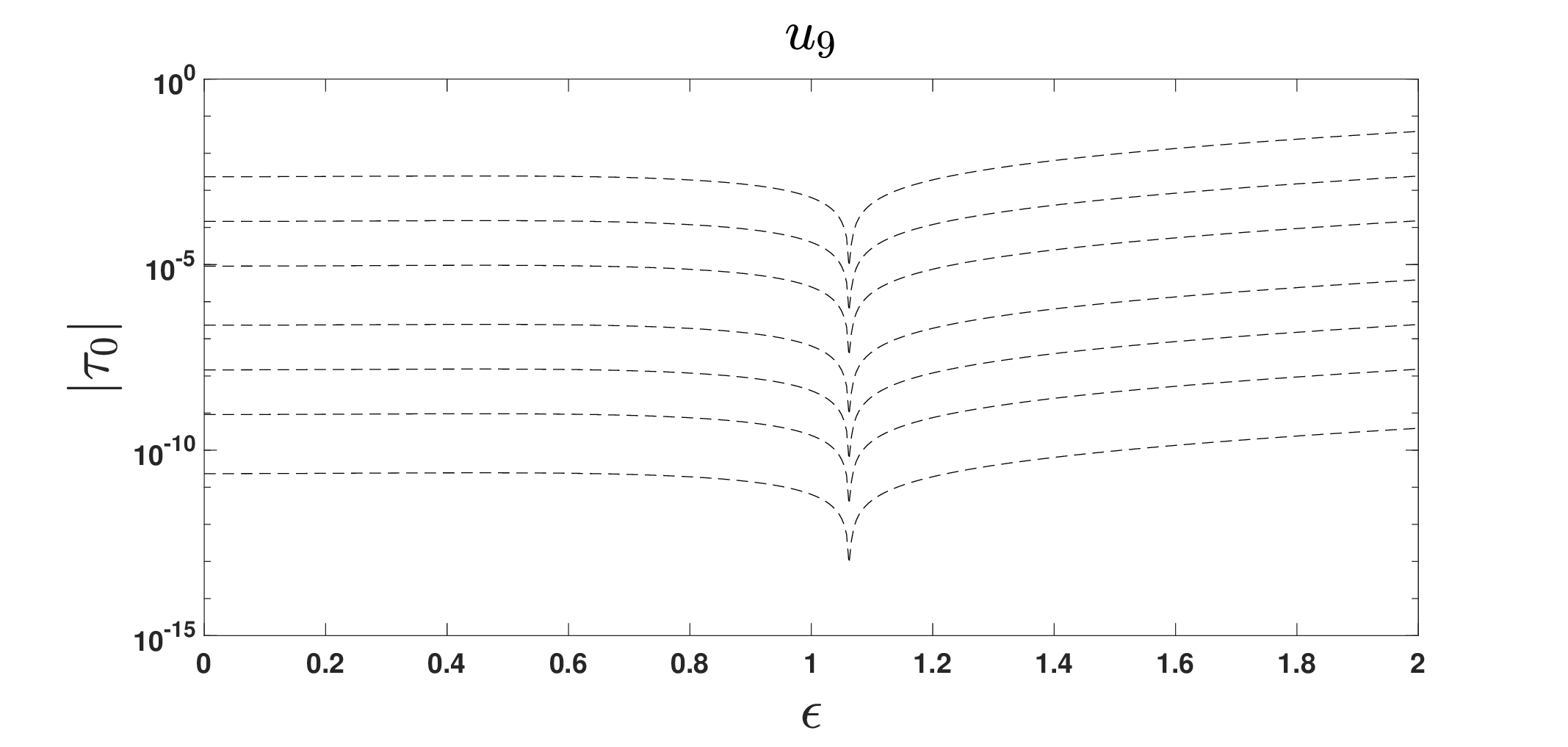}  
    \caption{Shape parameter ($\epsilon$) dependence on absolute LTEs ($|\tau_{0}|$) of RBF-HFD (order 4) formulas for 2D-Laplacian approximation of test functions (\ref{Test_function4_Laplacian}) and (\ref{Test_function9_Laplacian}) with step sizes, $ h = 0.2,\;0.1,\; 0.05,\; 0.02,\;0.01,\;0.005$\;and\;$ 0.002$, from top to bottom.} 
    \label{fig:Laplacian_Order4_LTE_EPS_1}
    \end{center}
\end{figure}
\subsection{Sixth order formula} 
\noindent The sixth order RBF-HFD formula for 2D-Laplacian is based on a stencil with local support sets (see \cite{collatznumerical, Satya_RBF_HFD})   
\begin{eqnarray}
S&=&\{  (x_0, y_0), (x_0+h,y_0), (x_0,y_0+h), (x_0-h,y_0), (x_0,y_0-h), (x_0+h, y_0+h), (x_0-h, y_0+h), \nonumber \\
&& (x_0-h, y_0-h), (x_0+h, y_0-h), (x_0+2h, y_0), (x_0,y_0+2h), (x_0-2h,y_0), (x_0,y_0-2h) \},\nonumber 
\\
\mu &=&\{ (x_0+h,y_0), (x_0,y_0+h), (x_0-h,y_0), (x_0,y_0-h), (x_0+h,y_0+h), (x_0-h,y_0+h), \nonumber \\ 
&& (x_0-h,y_0-h), (x_0+h,y_0-h) \}. \nonumber
\end{eqnarray}
In view of symmetry of stencil nodes layout and central difference approximation of 2D-Laplacian, the formula may be written as    
\begin{eqnarray}\label{RBF_HFD6_Scheme_Laplacian}
    \Delta u(x_0,y_0)&\thickapprox &  \alpha_{0}\; u_{0,0} + \alpha_{1} ( u_{1,0} +  u_{0,1} +  u_{-1,0} +  u_{0,-1} +u_{1,1} +  u_{-1,1} +  u_{-1,-1} +  u_{1,-1} ) \nonumber \\ 
  & & +\alpha_{2} ( u_{2,0} +  u_{0,2} +  u_{-2,0} +  u_{0,-2} ) + \beta_{1} \left( (\Delta u)_{1,0} +  (\Delta u)_{0,1}  +  (\Delta u)_{-1,0}  + (\Delta u)_{0,-1}) \right) \nonumber \\
  & & +\beta_{2} \left( (\Delta u)_{1,1} +  (\Delta u)_{-1,1}  +  (\Delta u)_{-1,-1}  + (\Delta u)_{1,-1}) \right).
\end{eqnarray}
The unknown weights in (\ref{RBF_HFD6_Scheme_Laplacian}) may be obtained by symbolically solving the corresponding linear system. To render the linear system amenable to symbolic computation in {\it Mathematica}, first we expand its non-constant elements using Taylor series expansion in powers of $\epsilon h$ up to order 24 and then symbolically solve the approximate linear system. Leading order expression for local truncation error ($\tau_0$) is obtained on substituting the weights and making use of Taylor series expansions of $u$ and $\Delta u$ values about reference node $(x_0,y_0)$. The analytical expressions of weights and LTEs are reported in Table \ref{Table3_RBF_HFD}. It is easy to observe that the weights of RBF-HFD (order 6) formulas derived for 2D-Laplacian converge to the weights of compact FD (order 6) \cite{collatznumerical} as $\epsilon \rightarrow 0$. The limiting weights are shown in the schematic stencil (\ref{Laplacian_RBF_HFD6_Weights_eps0}).

\begin{equation}\label{Laplacian_RBF_HFD6_Weights_eps0}
     \begin{tikzpicture}
        \node at (0,0) [rectangle, draw] (c1) {$\frac{-1}{46}$};
        \node at (1,0) [rectangle, draw] (c2) {$\frac{-5}{23}$};
        \node at (2,0) [rectangle, draw] (c3) {$\frac{-1}{46}$};
         \node at (-2,-1) [] (c25) { $\Delta u(x_0,y_0) \thickapprox \;\;\;\;\;$};
        \node at (0,-1) [rectangle,draw] (c4) {$\frac{-5}{23}$};
        \node at (2,-1) [rectangle,draw] (c6) {$\frac{-5}{23}$};
        \node at (3,-1) [] (c24) { $\Delta u +$};
        \node at (0,-2) [rectangle,draw] (c7) {$\frac{-1}{46}$};
        \node at (1,-2) [rectangle,draw] (c8) {$\frac{-5}{23}$};
        \node at (2,-2) [rectangle,draw] (c9) {$\frac{-1}{46}$};
        \draw (c1) -- (c9);
         \draw (c2) -- (c8);
         \draw (c3) -- (c7);
         \draw (c4) -- (c6);
         
          \node at (6,1) [rectangle, draw] (c18) {$\frac{9}{92}$};
          
         \node at (5,0) [rectangle, draw] (c15) {$\frac{12}{23}$};
        \node at (6,0) [rectangle, draw] (c16) {$\frac{12}{23}$};
        \node at (7,0) [rectangle, draw] (c17) {$\frac{12}{23}$};
         \node at (4,-1) [rectangle, draw] (c10) {$\frac{9}{92}$};
        \node at (5,-1) [rectangle, draw] (c11) {$\frac{12}{23}$};
        \node at (6,-1) [rectangle,draw] (c12) {$\frac{-105}{23}$};
        \node at (7,-1) [rectangle,draw] (c13) {$\frac{12}{23}$};
        \node at (8,-1) [rectangle,draw] (c14) {$\frac{9}{92}$};
        \node at (8.75,-1) [] (c26) {${\displaystyle\frac{u}{h^2}}$};
         \node at (5,-2) [rectangle,draw] (c19) {$\frac{12}{23}$};
        \node at (6,-2) [rectangle,draw] (c20) {$\frac{12}{23}$};
        \node at (7,-2) [rectangle,draw] (c21) {$\frac{12}{23}$};
         \node at (6,-3) [rectangle,draw] (c22) {$\frac{9}{92}$};
       \draw (c18) -- (c16);
       \draw (c15) -- (c12);
       \draw (c16) -- (c12);
       \draw (c17) -- (c12);
       \draw (c10) -- (c11);
       \draw (c11) -- (c12);
       \draw (c12) -- (c13);
       \draw (c13) -- (c14);
       \draw (c19) -- (c12);
       \draw (c20) -- (c12);
       \draw (c21) -- (c12);
       \draw (c22) -- (c20);
    \end{tikzpicture}
\end{equation}

\begin{figure}[h!]
    \begin{center}
    \includegraphics[width=48ex,height=40ex]{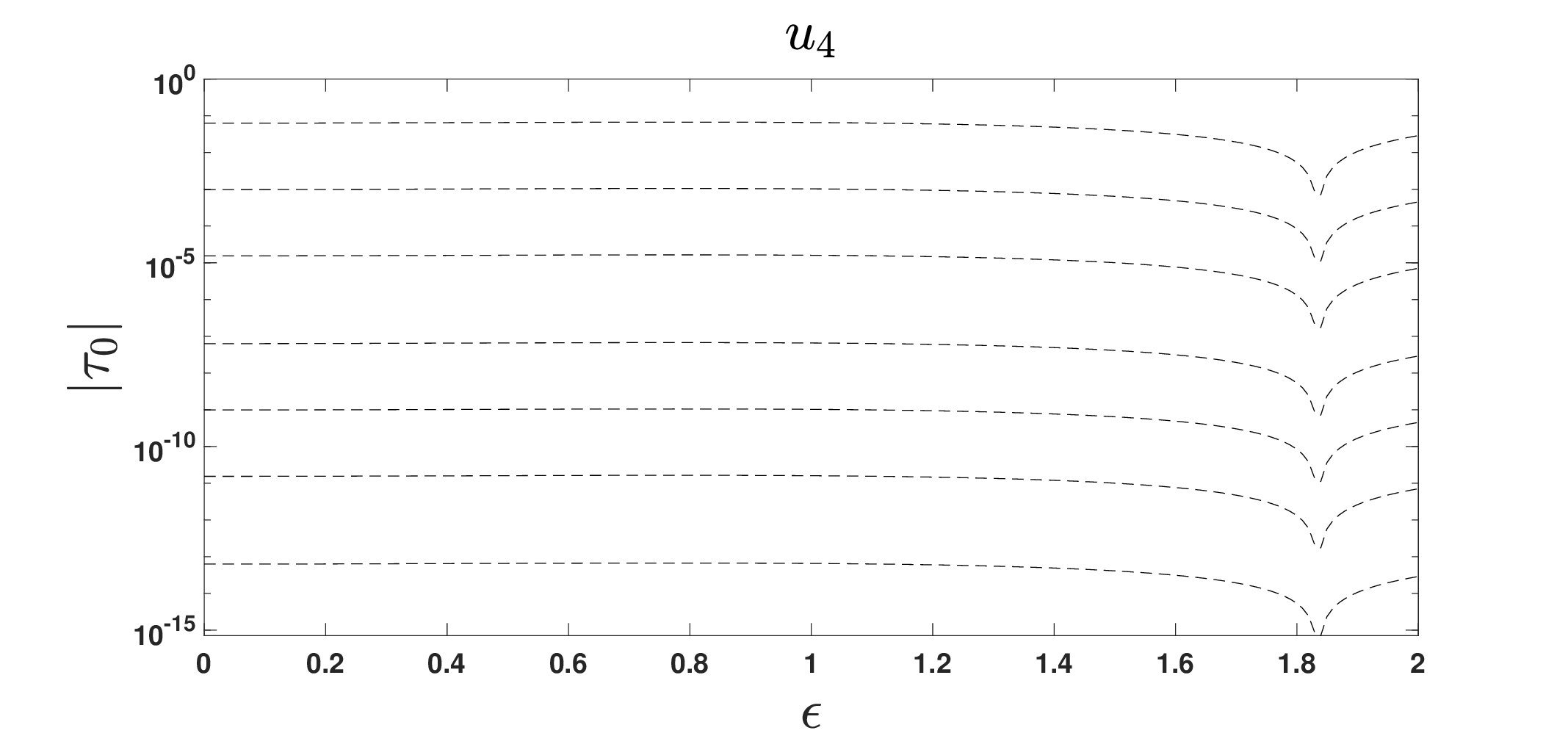}\includegraphics[width=48ex,height=40ex]{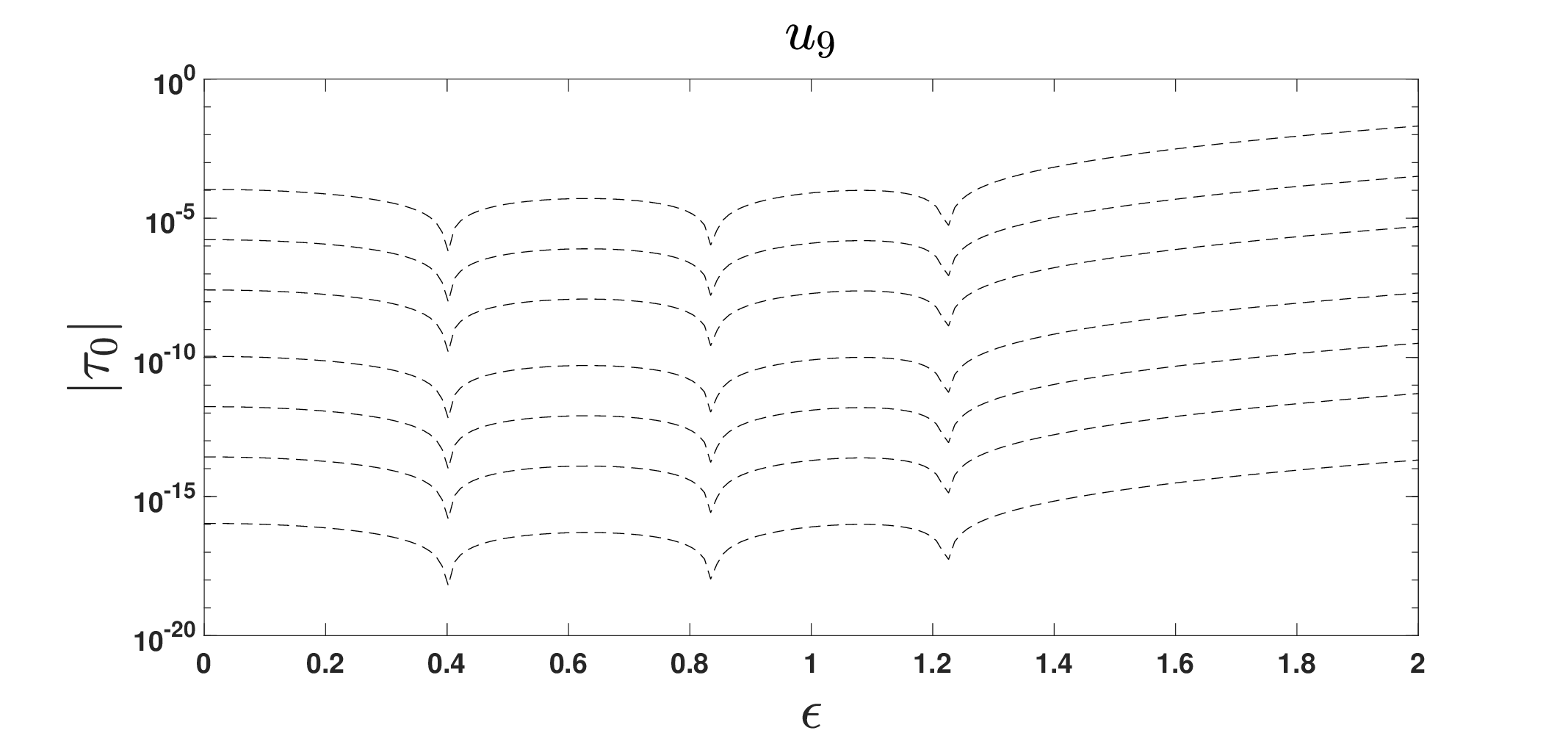}     
    \caption{Shape parameter ($\epsilon$) dependence on absolute LTEs ($|\tau_0|$) of RBF-HFD (order 6) formula for 2D-Laplacian approximation of test functions (\ref{Test_function4_Laplacian}) and (\ref{Test_function9_Laplacian}) with step sizes, $ h = 0.2,\;0.1,\; 0.05,\; 0.02,\;0.01,\;0.005$, $0.002$, from top to bottom.} 
    \label{fig:Laplacian_Order6_LTE_EPS_1}
    \end{center}
\end{figure}

\noindent Figures \ref{fig:Laplacian_Order4_LTE_EPS_1} and \ref{fig:Laplacian_Order6_LTE_EPS_1}, respectively, depict the dependence of shape parameter on absolute LTEs of fourth order (\ref{RBF_HFD4_Scheme_Laplacian}) and sixth order (\ref{RBF_HFD6_Scheme_Laplacian}) GA based formulas {\it w.r.t.} test functions $u_{4}$ and $u_{9}$ defined in (\ref{Test_function4_Laplacian}) and (\ref{Test_function9_Laplacian}). It is observed that for a chosen test function and reference point there exists a shape parameter value for which the absolute local truncation error is minimum and it is independent of step-size. In particular, the local truncation error for test function $u_{4}$ in Figure \ref{fig:Laplacian_Order4_LTE_EPS_1} and Figure 2 (Page. 115) of Wright and Fornberg \cite{wright2006scattered} has similar trend. However, our plot do not match exactly with the said figure because it includes the maximum error curves {\it w.r.t.} MQ-RBF whereas, we have shown the variation of local truncation error {\it w.r.t.} Gaussian-RBF for test function $u_{4}$. Thus, we may infer that even for fixed test function and reference point the optimal value of shape parameter may depend on the choice of radial basis function.  \\ 

\begin{figure}[h!]
    \begin{center}
\includegraphics[width=48ex,height=40ex]{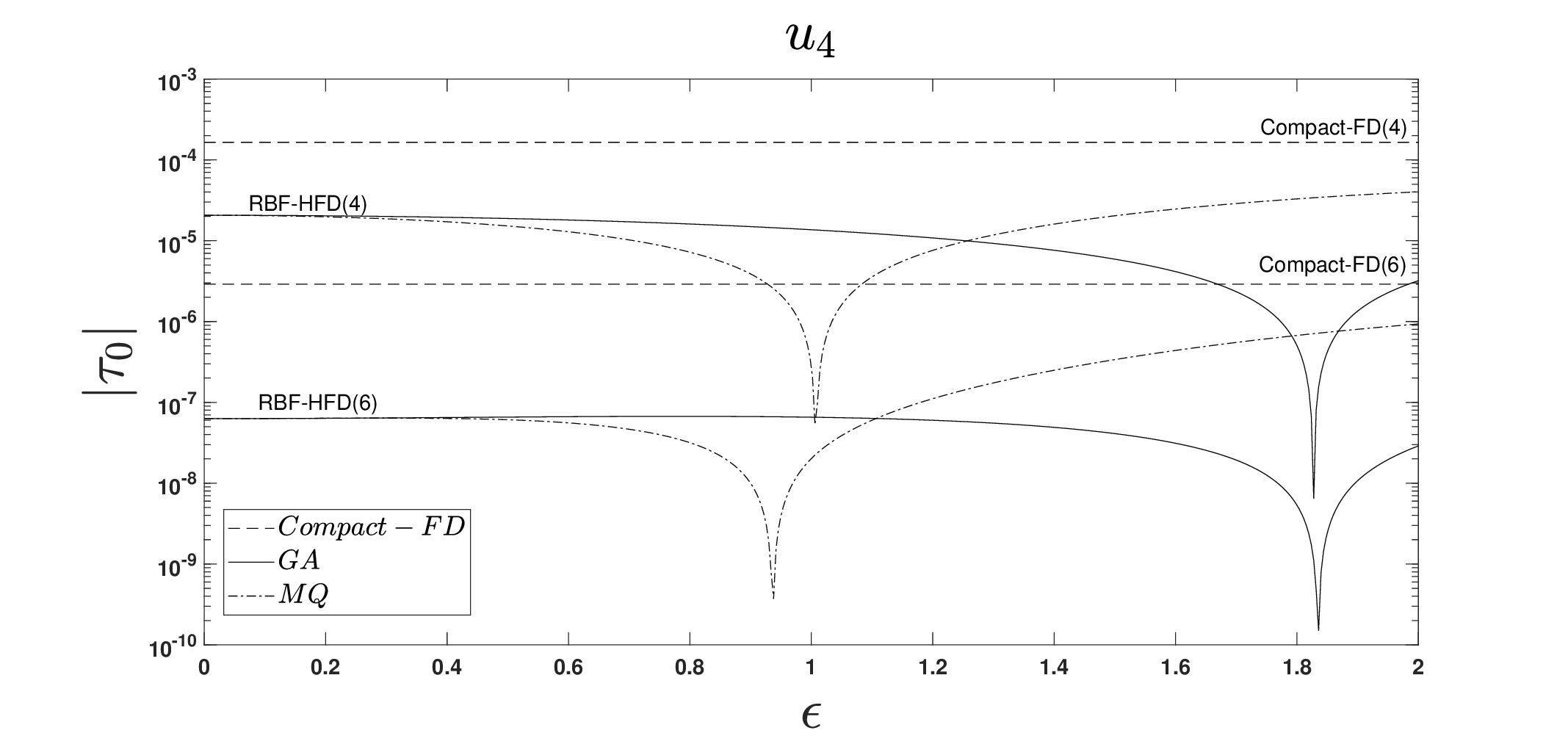}\includegraphics[width=48ex,height=40ex]{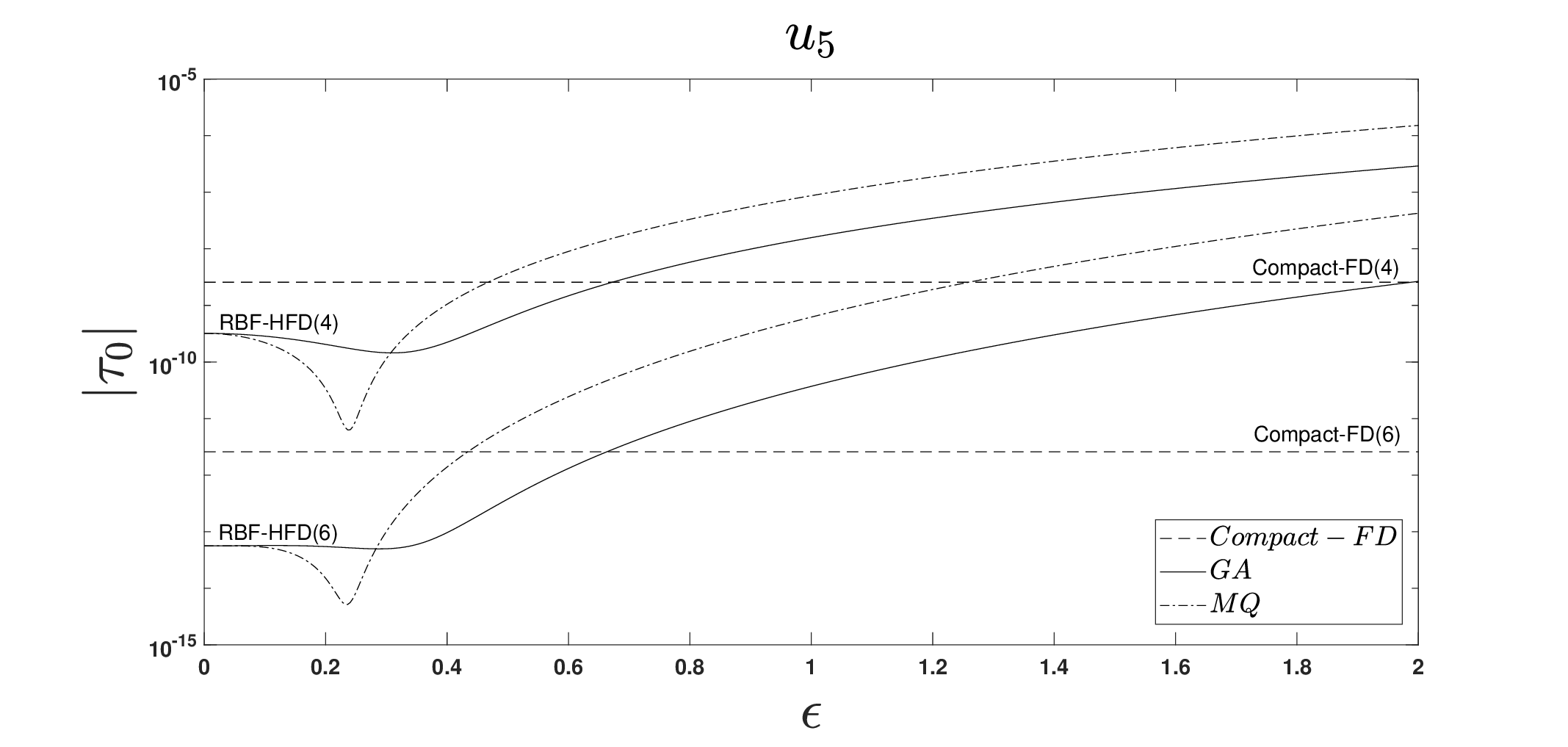}
     \includegraphics[width=48ex,height=40ex]{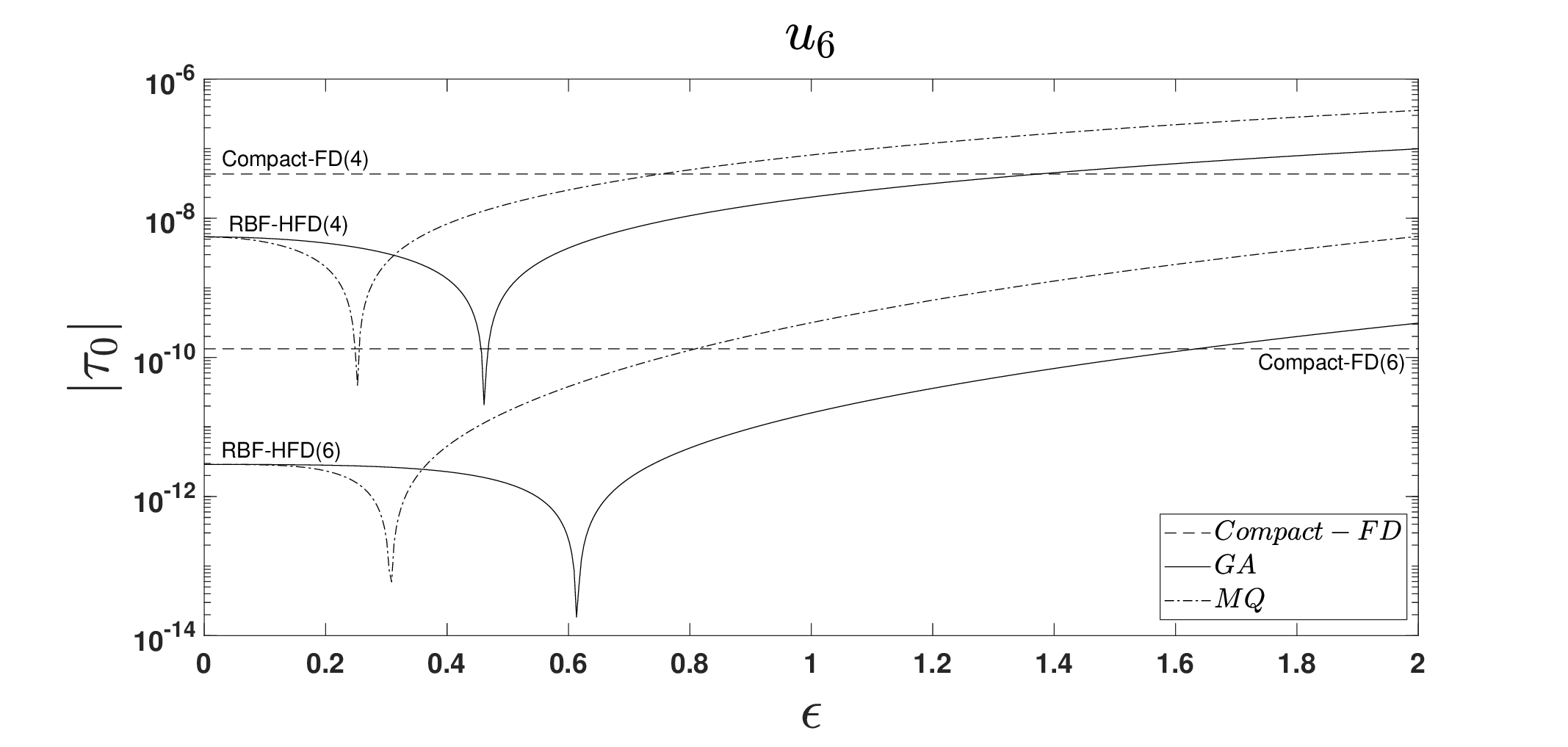}\includegraphics[width=48ex,height=40ex]{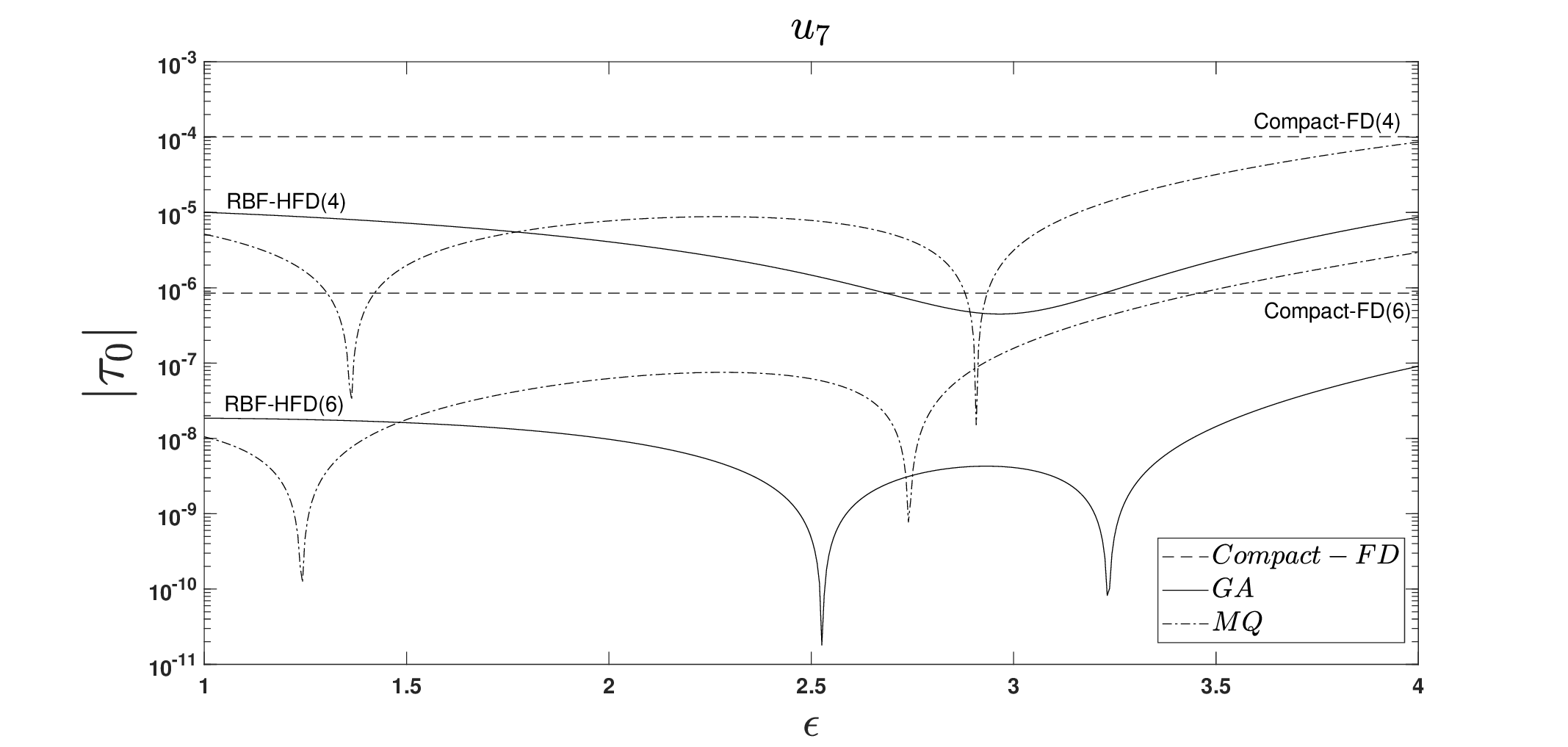}
\includegraphics[width=48ex,height=40ex]{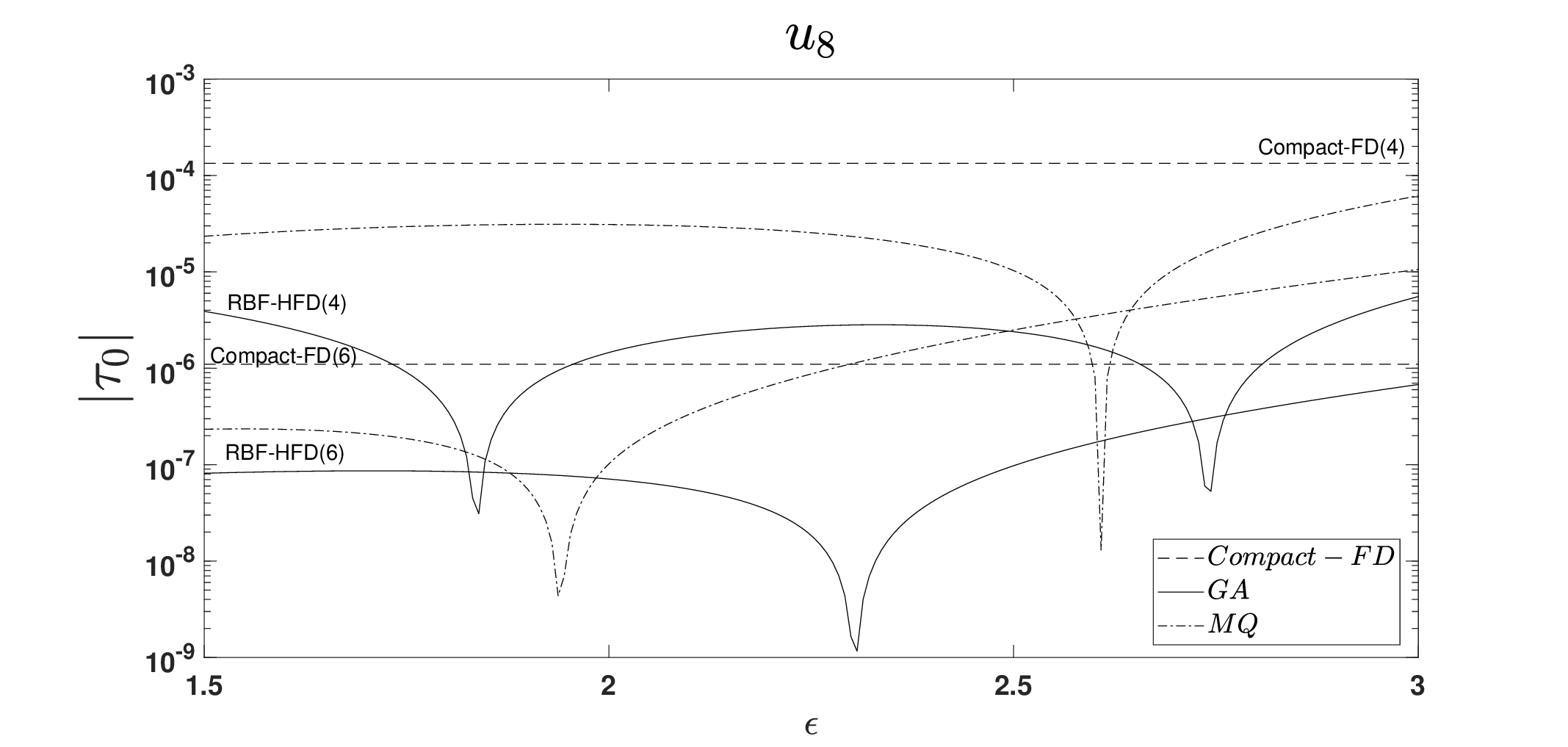}\includegraphics[width=48ex,height=40ex]{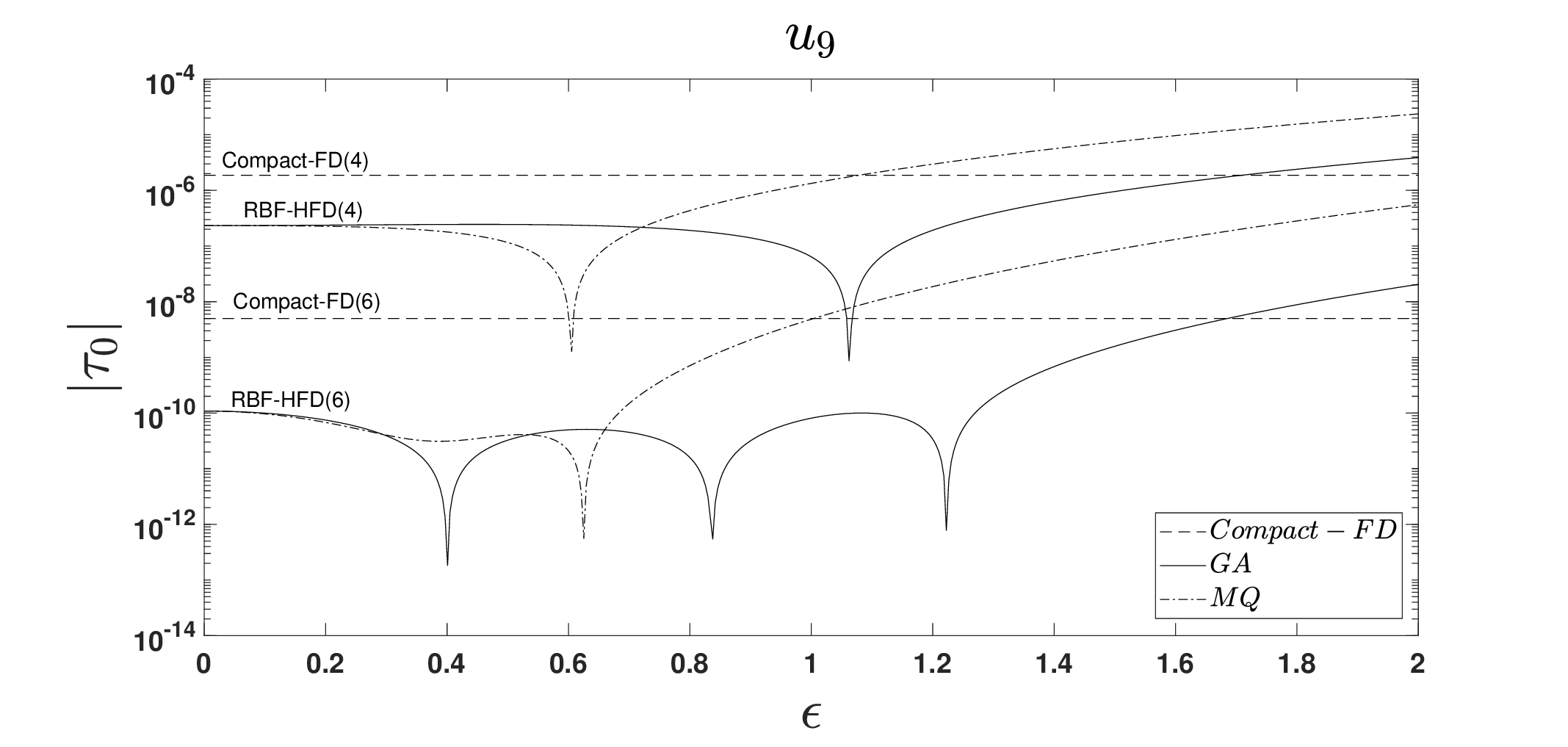}
\caption{Comparison of absolute LTEs ($\tau_{0}$) of compact FD (order 4 and 6), MQ and GA based RBF-HFD (order 4 and 6) formulas for 2D-Laplacian approximation of test functions (\ref{Test_function4_Laplacian})--(\ref{Test_function9_Laplacian}) with fixed step-size $h=0.02.$} 
    \label{fig:Laplacian_EPS_LTE_246_u3}
    \end{center}    
\end{figure}

\noindent We now compare the approximation accuracy of Gaussian based RBF-HFD (order 4 and 6) formulas for 2D-Laplacian with respective order compact FD schemes and MQ based RBF-HFD formulas \cite{Satya_RBF_HFD}. For this purpose we have chosen six test functions $u_4$, $u_5$, $u_6$, $u_7$, $u_8$, $u_9$ defined in equations (\ref{Test_function4_Laplacian})--(\ref{Test_function9_Laplacian}). Figure \ref{fig:Laplacian_EPS_LTE_246_u3}  
corresponds to the approximation accuracy of the test functions' 2D-Laplacian at specified reference points. In these figures, the horizontal dashed lines ($--$) represent the absolute LTEs of compact FD (order 4 and 6) schemes. The solid line curves represent absolute LTEs of the corresponding order Gaussian based RBF-HFD formulas, whereas dashed-dot ($-\cdot-$) line curves represent absolute LTEs of the corresponding order MQ based formulas. It is observed from these figures that, Gaussian based RBF-HFD formulas produce better accuracy in most of the cases as compared with compact FD formulas and MQ based formulas for certain ranges of shape parameter value.
\begin{landscape}
 \begin{table}[t]
    \caption{2D-Laplacian approximations using RBF-FD (order 2), RBF-HFD (order 4 and 6) formulas and comparison with classical FD2, Compact FD (order 4 and 6) schemes.}
\begin{center} 
\footnotesize{
\begin{tabular}{|c|c|c|} 
 \hline
 & & \\
   Stencil  & Schemes, weights and local truncation error (LTE) for RBF-HFD formulas & Weights and LTEs for  \\
& & compact FD formulas \\
\hline
& & \\
RBF-FD (Order 2) \cite{bayona2010rbf} (P.8287,Eq.17) & $(\Delta u)_{0,0} =\alpha_{0}\; u_{0,0} + \alpha_{1} ( u_{1,0} +  u_{0,1} +  u_{-1,0} +  u_{0,-1} ) $ & ( \cite{collatznumerical} Table 6, P. 542-543) \\
& &$\alpha_{0}=-\frac{4}{h^2},$ \\
 \begin{tikzpicture}
 \draw (0.75,-0.75) circle (0.15);  
\end{tikzpicture}& $\alpha_{0}=\frac{-4}{h^2}-{3\epsilon^2}+\frac{115\epsilon^6 h^4}{24}+O(\epsilon^8 h^6) $,\;$ \alpha_{1}=\alpha_{2} = \alpha_{3} =\alpha_{4} =\frac{1}{h^2}+\frac{3\epsilon^2}{4}-\frac{115\epsilon^6 h^4}{96}+O(\epsilon^8 h^6) $&$\alpha_{1}=\alpha_{2} = \alpha_{3} =\alpha_{4} =\frac{1}{h^2},$\\
\begin{tikzpicture}
  \draw (0,0) circle (0.15);
  \draw[line width=1mm] (0.75,0) circle (0.15);  
  \draw (1.5,0) circle (0.15);
   \end{tikzpicture} & $ \tau_{0}=\frac{h^2}{12} \left(u^{(4,0)}(\mathbf{x})+u^{(0,4)}(\mathbf{x})\right)+{9 \epsilon^2}(u^{(2,0)}(\mathbf{x})+u^{(0,2)}(\mathbf{x})),\;\;\text{where $\mathbf{x}=(x_0,y_0)$}$ &  $ \tau_{0} \approx -\dfrac{h^2}{12} \left(u^{(4,0)}(\mathbf{x})+u^{(0,4)}(\mathbf{x})\right)$\\
 \begin{tikzpicture}
     \draw (0.75,-0.75) circle (0.15);
 \end{tikzpicture}& & \\
\hline
& &\\
 RBF-HFD (Order 4) &$(\Delta u)_{0,0}=\alpha_{0}\; u_{0,0} + \alpha_{1} ( u_{1,0} +  u_{0,1} +  u_{-1,0} +  u_{0,-1} ) + \alpha_{2} ( u_{1,1} +  u_{-1,1} +  u_{-1,-1} +  u_{1,-1} ) $ & ( \cite{collatznumerical} Table 6, P. 542-543)\\ 
 & $+ \beta_{1} \left( (\Delta u)_{1,0} +  (\Delta u)_{0,1}  +  (\Delta u)_{-1,0}  + (\Delta u)_{0,-1}) \right)$ & \\
& $ \alpha_{0} = -\frac{5}{h^2} - \frac{19 \epsilon^2}{24} + \frac{4069 \epsilon^4 h^2}{1440} + \frac{9791 \epsilon^6 h^4}{8640}+O(\epsilon^8 h^6)$ & $ \alpha_{0} = -\frac{5}{h^2},$ \\
\begin{tikzpicture}
\draw (0,-0.75) circle (0.15);
  \draw (0.75,-0.75) circle (0.15); 
   \draw (0.75,-0.75) circle (0.25); 
  \draw (1.5,-0.75) circle (0.15); 
  \end{tikzpicture} &$\alpha_{1}=\frac{1}{h^2}-\frac{13 \epsilon^2}{24}-\frac{2273\epsilon^4 h^2}{1440}-\frac{4459 \epsilon^6 h^4}{8640} +O(\epsilon^8 h^6) $ & $\alpha_{1}=\frac{1}{h^2},$\\
\begin{tikzpicture}
  \draw (0,0) circle (0.15);
   \draw (0,0) circle (0.25);
  \draw[line width=1mm] (0.75,0) circle (0.15);  
  \draw (1.5,0) circle (0.15);
  \draw (1.5,0) circle (0.25);
\end{tikzpicture} &$\alpha_{2}=\frac{1}{4 h^2}+\frac{71\epsilon^2}{96}+\frac{5023\epsilon^4 h^2}{5760}+\frac{1609\epsilon^6 h^4}{6912}+O(\epsilon^8 h^6) $ & $\alpha_{2}=\frac{1}{4 h^2},$ \\
 \begin{tikzpicture}
  \draw (0,0.75) circle (0.15);
  \draw (0.75,0.75) circle (0.15);
  \draw (0.75,0.75) circle (0.25);
  \draw (1.5,0.75) circle (0.15);
 \end{tikzpicture}&$\beta_{1}=-\frac{1}{8}-\frac{15 \epsilon^2 h^2}{64}-\frac{213 \epsilon^4 h^4}{1280}+\frac{373 \epsilon^6 h^6}{7680}+O(\epsilon^8 h^8) $& $\beta_{1}=-\frac{1}{8},$ \\
 & $ \tau_{0}= \frac{h^4}{480}( -240 \epsilon^4 (u^{(2,0)}(\textbf{x})+u^{(0,2)}(\textbf{x}))-75 \epsilon^2( u^{(4,0)}(\textbf{x})+ u^{(0,4)}(\textbf{x}))$&\\
 &$-3(u^{(6,0)}(\mathbf{x})+u^{(0,6)}(\mathbf{x}))+5(u^{(2,4)}(\mathbf{x})+ u^{(4,2)}(\mathbf{x}))+130 \epsilon^2 u^{(2,2)}(\mathbf{x})) + O(h^6 P_{3}(\epsilon^2))\;\;\text{where $\mathbf{x}=(x_0,y_0)$}$ & $\tau_{0} \approx \dfrac{h^4}{60}\left(3((u^{(6,0)}(\mathbf{x})+u^{(0,6)}(\mathbf{x})\right))$ \\
& &$-5(u^{(2,4)}(\mathbf{x})+u^{(4,2)}(\mathbf{x})))$\\
 & &\\ 
\hline 
& & (\cite{collatznumerical} Table 6, P. 542-543) \\
RBF-HFD (Order 6)&$(\Delta u)_{0,0}=\alpha_{0}\; u_{0,0} + \alpha_{1} ( u_{1,0} +  u_{0,1} +  u_{-1,0} +  u_{0,-1} +u_{1,1} +  u_{-1,1} +  u_{-1,-1} +  u_{1,-1} ) +\alpha_{2} ( u_{2,0} +  u_{0,2} +  u_{-2,0} +  u_{0,-2} )$ & \\
&$ + \beta_{1} \left( (\Delta u)_{1,0} +  (\Delta u)_{0,1}  +  (\Delta u)_{-1,0}  + (\Delta u)_{0,-1}) \right)+\beta_{2} \left( (\Delta u)_{1,1} +  (\Delta u)_{-1,1}  +  (\Delta u)_{-1,-1}  + (\Delta u)_{1,-1}) \right)$ &$\alpha_{0}=-\frac{105}{23 h^2},$ \\
\begin{tikzpicture}
      \draw (0.75,-1.5) circle (0.15); 
\end{tikzpicture} &$ \alpha_{0}=\frac{-105}{23 h^2} +\frac{61761 \epsilon^2}{83582}-\frac{2668325501\epsilon^4 h^2}{1518684940}-\frac{2725280260345847\epsilon^6 h^4}{289742306277900}+O(\epsilon^8 h^6)$ &$\alpha_{1}=\frac{12}{23 h^2}$\\
&$\alpha_{1}=\frac{12}{23 h^2}-\frac{17646 \epsilon^2}{208955}-\frac{3697010966 \epsilon^4 h^2}{13288493225}+\frac{1329393861667 \epsilon^6 h^4}{724355765694750}
+O(\epsilon^8 h^6) $&$ \alpha_{2}=\frac{9}{92 h^2}$ \\
\begin{tikzpicture}
\draw (0,-0.75) circle (0.15);
\draw (0,-0.75) circle (0.25);
  \draw (0.75,-0.75) circle (0.15); 
   \draw (0.75,-0.75) circle (0.25); 
  \draw (1.5,-0.75) circle (0.15); 
   \draw (1.5,-0.75) circle (0.25); 
  \end{tikzpicture} & $ \alpha_{2}=\frac{9}{92 h^2}-\frac{26469 \epsilon^2}{1671640}-\frac{24912958377 \epsilon^4 h^2}{212615891600}+\frac{4535043666647521 \epsilon^6 h^4}{1931615375186000}+O(\epsilon^8 h^6)$&$\beta_{1}=-\frac{5}{23},$\\
\begin{tikzpicture}
 \draw (-1,0) circle (0.15);
  \draw (0,0) circle (0.15);
   \draw (0,0) circle (0.25);
  \draw[line width=1mm] (0.75,0) circle (0.15);  
  \draw (1.5,0) circle (0.15);
  \draw (1.5,0) circle (0.25);
   \draw (2.5,0) circle (0.15);
\end{tikzpicture}&$ \beta_{1}=\frac{-5}{23}+\frac{97053 \epsilon^2 h^2}{835820} +\frac{19096725789 \epsilon^4 h^4}{106307945800}-\frac{2108622124065179 \epsilon^6 h^6}{724355765694750}
+O(\epsilon^8 h^8)      $&$\beta_{2}=-\frac{1}{46}$\\
\begin{tikzpicture}
\draw (0,0.75) circle (0.15);
\draw (0,0.75) circle (0.25);
  \draw (0.75,0.75) circle (0.15); 
   \draw (0.75,0.75) circle (0.25); 
  \draw (1.5,0.75) circle (0.15); 
   \draw (1.5,0.75) circle (0.25); 
    \end{tikzpicture}& $\beta_{2}=\frac{-1}{46}-\frac{61761 \epsilon^2 h^2}{1671640}-\frac{8234946399 \epsilon^4 h^4}{30373698800}-\frac{828029839457483\epsilon^6 h^6}{5794846125558000}+O(\epsilon^8 h^8)$&\\
 & $ \tau_{0}= \frac{h^6}{3052560} (-8605128 \epsilon^6 (u^{(2,0)}(\mathbf{x})+u^{(0,2)}(\mathbf{x}))-1371468 \epsilon^4( u^{(4,0)}(\mathbf{x})+ u^{(0,4)}(\mathbf{x}))+1817(u^{(8,0)}(\mathbf{x}))+u^{(0,8)}(\mathbf{x}))$&$\tau_{0} \approx \dfrac{h^6}{1260}\left(-23((u^{(8,0)}(\mathbf{x})+u^{(0,8)}(\mathbf{x}))\right)+$\\
\begin{tikzpicture}
     \draw (0.75,1.5) circle (0.15); 
\end{tikzpicture}&$-123522\epsilon^2(u^{(2,4)}(\mathbf{x})+ u^{(4,2)}(\mathbf{x}))-3318(u^{(2,6)}(\mathbf{x})+u^{(6,2)}(\mathbf{x}))-13644 \epsilon^4 u^{(2,2)}(\mathbf{x}))  + O(h^8 P_{4}(\epsilon^2)),$\text{where $\mathbf{x}=(x_0,y_0)$} & \ \ \ \ \ \ \ \ \ \ \ \ \  \ \  $\left(42(u^{(2,6)}(\mathbf{x})+ u^{(6,2)}(\mathbf{x})\right)$\\
 & &\\

 & &\\
 \hline
\end{tabular}
}
\label{Table3_RBF_HFD}
\end{center}
\end{table}
\end{landscape}
\section{Optimal shape parameter}\label{Optimal_shape_parameter}
  
\noindent Following Satyanarayana {\it et al.} \cite{Satya_RBF_HFD} we make use of an optimization algorithm, proposed by Bayona {\it et al.} \cite{bayona2010rbf} for RBF-FD approximations, for computing optimal values of shape parameter associated to GA based RBF-HFD formulas derived in Sections \ref{First_derivative_approximation}, \ref{Second_derivative_approximation}, and \ref{Two_dimensional_Laplacian_approximation}. Local truncation error at reference point $x_{0}$ for an RBF-HFD (order $m$) formula may be written in the form 
\begin{eqnarray}\label{Opt_eps}
          \tau_{0}=h^{m} P_{n}(z) + O(h^{m+2} P_{n+1}(z)),\;\mbox{where} \;  z=\epsilon^{2}\; 
          \mbox{and}\;   P_{n}(z) = \sum_{i=0}^{n} a_{i} z^{i},\; a_{i} \in \mathbb{R}.
\end{eqnarray}

\noindent The coefficients $a_{i}$ in polynomial $P_n(z)$ depends on the choice of test function and the reference point. Then the algorithm requires us to find some $z_c \in \mathbb{R}^{+} $ such that either $P_n(z_c)$ is zero, or $z_c$ is the global minima of $P_n(z)$. Then the optimal shape parameter value is $\epsilon^{*} = \sqrt{z_c}$. \\

\noindent Sixth order RBF-HFD formula (\ref{1D_Six_Order_Scheme}) for approximating first derivative leads to local truncation error
\begin{eqnarray}
          \tau_{0}=h^{6} P_{3}(z) + O(h^{8} P_{4}(z)),~~ P_{3}(z) =a_{0} + a_{1}z + a_{2}z^{2} + a_{3}z^{3}, \label{Opt_LTE_6}
\end{eqnarray}
\noindent where $a_{0}=\frac{1}{120} u^{(7)}(x_0)$, $a_{1}=\frac{42}{120} u^{(5)}(x_0)$, $a_{2}=\frac{420}{120} u^{(3)}(x_0)$ and $a_{3}=\frac{840}{120} u'(x_0)$. Then for the test function $u_{1}$ defined in (\ref{u1}) with reference point $x_{0}= 0.4$, the coefficients assume the values $$a_{0} = \frac{16783616k_{1}}{78125}+\frac{2056992k_{2}}{3125}, 
     a_{1} = \frac{10752k_{2}}{25}-\frac{6256992k_{1}}{3125}, a_{2} = -\frac{5376k_{1}}{25}-2016k_{2}, \; a_{3} = 672k_{1},$$ where $k_{1}=\cos(4/25)$ and $k_{2}=\sin(4/25)$. Solving the cubic equation $P_{3}(z) = 0$, we obtain two positive roots and one negative root,  $z_{1} =2.077764948 $, $z_{2} =0.1603514855 $ and $z_{3} = -1.433978051$. Thus the optimal shape parameter is $\epsilon^{*}=\sqrt{z_{2}}=0.4004391108$.  \\
     
     \noindent Second derivative approximation using RBF-HFD (order 6) formula (\ref{RBF_HFD6_2D}) leads to local truncation error of the form (\ref{Opt_LTE_6}). Then for test function $u_{1}$ at reference point $x_{0}= 0.4$, the cubic polynomial coefficients assume values
$a_{0} = \frac{33370624k_{1}}{15625}+\frac{521915536k_{2}}{390625}, \;
     a_{1} = \frac{96978816k_{2}}{15625}-\frac{727056k_{1}}{125}, \;
      a_{2} = -\frac{145152k_{1}}{25}-\frac{5476464k_{2}}{625} \;\mbox{and}\;
       a_{3} = 5120k_{1} -\frac{8192k_{2}}{5}$, where $k_1= \cos(4/25)$ and $k_2=\sin(4/25)$.
Two roots of the polynomial equation $P_3(z)=0$ are positive and third is negative, $z_{1} = 1.877718044 $,\; $z_{2} = 0.3487463007 $ and $z_{3} =-0.73949091922  $. Thus the optimal shape parameter is  $\epsilon^{*}=\sqrt{z_{2}}=0.5905474585 $. \\ 
\begin{table}[h!]
   \begin{center} 
\caption{Optimal shape parameters ($\epsilon^{*}$) for GA based RBF-HFD formulas to approximate the first and second derivatives of test functions $u_{1}$ and $u_{2}$.} 
\begin{tabular}{|c|c|c|c|c|c|} 
\hline
& Test function  & $\epsilon^{*}$  & $\epsilon^{*}$  & $\epsilon^{*}$  & $\epsilon^{*}$ \\
&    &   RBF-HFD(4) &   RBF-HFD(6) &   RBF-HFD(8) &   RBF-HFD(10) \\
\hline
First
& $u_{1}$ at  $x_0= 0.4$ & 1.13112 & 0.39038 & 0.75074 & 0.37036 \\
 Derivative & $u_{2}$ at  $x_0=0.25$ & 0.87086 & 0.57056 & 0.65064 & 0.35034 \\
\hline
Second 
&  $u_{1}$ at  $x_0=0.4$ & 1.13112 & 0.59058 & 0.85084 & 0.49048 \\
 Derivative & $u_{2}$ at  $x_0=0.25$ & 0.65064 & 0.55054 & 0.45044 & 0.33032\\
\hline
\end{tabular}\label{Table4_RBF_HFD_Opt}  
\end{center}
\end{table}  \\   
\noindent The optimal shape parameter for RBF-HFD (order 4) formula (\ref{RBF_HFD4_Scheme_Laplacian})  to approximate 2D-Laplacian is computed using its local truncation error expression. With respect to test function $u_{4}$ at reference point $(x_{0},y_{0})=(0.25,0.25)$, the leading order term is a quadratic polynomial 
\begin{equation}
P_{2}(z) =a_{0} + a_{1}\;z + a_{2}\;z^{2} 
       \label{4thoptLap}
\end{equation}
\noindent where $ a_{0} = 99 \pi^{5} k +\dfrac{985\pi^{3} k}{2} + \dfrac{5129\pi k}{16}$, 
     $a_{1} =-{470 \pi^{3}k}- \dfrac{1205 \pi k}{2}$, $a_{2} = 240 \pi k $ and $k=\sqrt{2} e^{(-1/16)}$. The roots of $P_2(z)=0$ are $z_{1} =18.49934524$ and $z_{2} = 3.339046708$. Then the optimal value of shape parameter is  $\epsilon^{*}=\sqrt{z_{2}}=1.827305861$.  From equation (\ref{RBF_HFD4_Scheme_Laplacian}) the optimal shape parameter with respect to test function $u_{5}$ at reference point $(x_{0},y_{0})=(0,0)$ is obtained from the roots of quadratic polynomial $P_2(z)$ whose coefficients are
$$ a_{0} = \frac{1756875000}{30664297}, \;
     a_{1} = - \frac{33033466796875}{3004150512793}, \;
     a_{2} = \frac{282820141259765625}{294313621587817417}.$$ Then roots of  equation $P_2(z)=0$ are complex valued. Therefore, we look for roots of derivative of $P_2(z)$
\begin{equation}
    2 a_{2} z + a_{1}=0. \label{4thoptLapder2}
\end{equation}
\noindent Then root of equation \eqref{4thoptLapder2} is $z = -a_1/2a_2 = 0.09596096811$. Thus,  $\epsilon^{*} = \sqrt{z} = 0.3097756739$ is the optimal shape parameter value. In Table 5, we report the optimal shape parameter values for GA based RBF-HFD (order 4, 6) formulas for 2D-Laplacian associated to test functions $u_{4}$ to $u_{9}$ defined in (\ref{Test_function4_Laplacian})--(\ref{Test_function9_Laplacian}). 
\begin{table}[h!]
   \begin{center} 
\caption{Optimal shape parameter values ($\epsilon^{*}$) for GA based RBF-HFD formulas for Laplacian of test functions defined in (\ref{Test_function4_Laplacian})--(\ref{Test_function9_Laplacian})} 
\begin{tabular}{|c|c|c|c|} 
\hline
& Test function   & $\epsilon^{*}$  & $\epsilon^{*}$   \\
&    &   RBF-HFD(4) &   RBF-HFD(6) \\
\hline
& $u_{4}$ at $(x_0,y_0)=(0.25,0.25)$ & 1.8277 & 1.8392  \\
& $u_{5}$ at $(x_0,y_0)=(0,0)$ & 0.2965 &  0.2914  \\
Laplacian 
& $u_{6}$ at $(x_0,y_0)=(0.1,0.2)$ & 0.4609 & 0.6130  \\
& $u_{7}$ at $(x_0,y_0)=(0.1,0.2)$ & 2.9748 & 2.5226  \\
& $u_{8}$ at $(x_0,y_0)=(0.1,0.2)$ & 1.8356 & 2.3065  \\
& $u_{9}$ at $(x_0,y_0)=(0.1,0.2)$ & 1.0621 & 0.4020  \\
 \hline
\end{tabular} 
\label{Table5_RBF_HFD_Opt} 
\end{center}
\end{table}
\section{Conclusions}\label{Conclusions}
\noindent In this article, analytical expressions of weights and local truncation errors for Gaussian based RBF-HFD formulas for first derivative, second derivative and 2D-Laplacian are derived. Analytical expressions help remove instabilities and ill-conditioning arising from direct numerical computation. The convergence properties of the formulas are validated for standard test functions. The main conclusions of the study are as follows :   
\begin{itemize}
     \item Compact FD formulas may be obtained as flat limit of the Gaussian based RBF-HFD formulas.
     \item  For some choices of test functions, Gaussian based formulas yield better approximation accuracy as compared to compact FD schemes and MQ based formulas.
    \item Optimal shape parameter values for each formula is computed for several test functions and step-sizes.
\end{itemize}
\bibliographystyle{elsarticle-num}
\bibliography{biblatex}

\begin{thebibliography}{10}
\expandafter\ifx\csname url\endcsname\relax
  \def\url#1{\texttt{#1}}\fi
\expandafter\ifx\csname urlprefix\endcsname\relax\def\urlprefix{URL }\fi
\expandafter\ifx\csname href\endcsname\relax
  \def\href#1#2{#2} \def\path#1{#1}\fi

\bibitem{wright2006scattered}
G.~B. Wright, B.~Fornberg, Scattered node compact finite difference-type formulas generated from radial basis functions, J. Comput. Phys. 212 (2006) 99--123.

\bibitem{collatznumerical}
L.~Collatz, The numerical treatment of differential equations. \uppercase{B}erlin: Springer; (1960).

\bibitem{lele1992compact}
S.~K. Lele, Compact finite difference schemes with spectral-like resolution, J. Comput. Phys. 103 (1992) 16--42.

\bibitem{fornberg2004stable}
B.~Fornberg, G.~B. Wright, Stable computation of multiquadric interpolants for all values of the shape parameter, Comput. Math. Appl. 48 (2004) 853--867.

\bibitem{fornberg2011stable}
B.~Fornberg, E.~Larsson, N.~Flyer, Stable computations with {G}aussian radial basis functions, SIAM J. Sci. Comput. 33 (2011) 869--892.

\bibitem{fornberg2013stable}
B.~Fornberg, E.~Lehto, C.~Powell, Stable calculation of {G}aussian-based \uppercase{rbf-fd} stencils, Comput. Math. Appl. 65 (2013) 627--637.

\bibitem{wright2017stable}
G.~B. Wright, B.~Fornberg, Stable computations with flat radial basis functions using vector-valued rational approximations, J. Comput. Phys. 331 (2017) 137--156.

\bibitem{bayona2010rbf}
V.~Bayona, M.~Moscoso, M.~Carretero, M.~Kindelan, \uppercase{RBF-FD} formulas and convergence properties, J. Comput. Phys. 229 (2010) 8281--8295.

\bibitem{bayona2011optimal}
V.~Bayona, M.~Moscoso, M.~Kindelan, Optimal constant shape parameter for multiquadric based \uppercase{RBF-FD} method, J. Comput. Phys. 230 (2011) 7384--7399.

\bibitem{bayona2012optimal}
V.~Bayona, M.~Moscoso, M.~Kindelan, Optimal variable shape parameter for multiquadric based \uppercase{RBF-FD} method, J. Comput. Phys. 231 (2012) 2466--2481.

\bibitem{bayona2012gaussian}
V.~Bayona, M.~Moscoso, M.~Kindelan, Gaussian \uppercase{RBF-FD} weights and its corresponding local truncation errors, Eng. Anal. Bound. Elem. 36 (2012) 1361--1369.

\bibitem{Satya_RBF_HFD}
C.~Satyanarayana, M.~K. Yadav, M.~Nath, Multiquadric based {RBF-HFD} approximation formulas and convergence properties, Eng. Anal. Bound. Elem. 160 (2024) 234--257.

\bibitem{song2024computing}
Y.~Song, M.~Barfeie, F.~Soleymani, Computing compact finite difference formulas under radial basis functions with enhanced applicability, Appl. Numer. Math 201 (2024) 370--386.

\bibitem{ding2005error}
H.~Ding, C.~Shu, D.~Tang, Error estimates of local multiquadric-based differential quadrature \uppercase{(LMQDQ)} method through numerical experiments, Int. J. Numer. Methods Eng. 63 (2005) 1513--1529.

\bibitem{chandhini2007local}
G.~Chandhini, Y.~Sanyasiraju, Local \uppercase{rbf-fd} solutions for steady convection--diffusion problems, Int. J. Numer. Methods Eng. 72 (2007) 352--378.

\bibitem{sanyasiraju2008local}
Y.~Sanyasiraju, G.~Chandhini, Local radial basis function based gridfree scheme for unsteady incompressible viscous flows, J. Comput. Phys. 227 (2008) 8922--8948.

\end{thebibliography}
\end{document}